%% file: qbxpaper.tex
\newif\ifelsart

\newif\ifsubmit
\submittrue

\ifelsart
\documentclass[a4paper,11pt,english]{elsarticle}
\else
\documentclass[a4paper,11pt,english]{article}
\fi
\usepackage{geometry}
\usepackage[T1]{fontenc}
\usepackage{babel}
\usepackage{graphicx}
\usepackage{caption,subcaption}
\usepackage{bm}
\usepackage{amsmath,amssymb}
\usepackage{jabbrv,doi}
\usepackage[title]{appendix}
\usepackage{hyperref}
\hypersetup{
  pdfborder = 0 0 0,
  colorlinks=false,
  linkcolor=black,
  citecolor=black,
  filecolor=black,
  urlcolor=black,
} %

\ifsubmit
\usepackage{tikzexternal,tikzscale}
\graphicspath{{ext}}
\else
\usepackage{tikz,tikzscale,pgfplots}
\usetikzlibrary{plotmarks}
\usetikzlibrary{external}
\pgfplotsset{compat=newest} 
\pgfplotsset{plot coordinates/math parser=false} 
\tikzset{mark size=1/2}
\tikzset{font=\footnotesize}
\fi

\AtEndPreamble{
  \LetLtxMacro{\oldincludegraphics}{\includegraphics}
  \newcommand{\myincludegraphics}[2][]{\tikzsetnextfilename{#2}\oldincludegraphics[#1]{#2}}
  \LetLtxMacro{\includegraphics}{\myincludegraphics}
}

\input{defs.tex}

\input{acronyms.tex}

\newcommand{\thetitle}{A fast integral equation method for solid
  particles in viscous flow using quadrature by expansion}

\ifelsart
\begin{document}
\begin{frontmatter}
\author[kth]{Ludvig af Klinteberg\corref{cor1}}
\ead{ludvigak@kth.se}
\cortext[cor1]{Corresponding author}
\author[kth]{Anna-Karin Tornberg}
\address[kth]{
  Numerical Analysis, Department of Mathematics,\\
  KTH Royal Institute of Technology, 100 44 Stockholm, Sweden
}
\title{\thetitle}

\else
\usepackage{authblk}
\author{Ludvig af Klinteberg\thanks{Corresponding author. Email:
    \href{ludvigak@kth.se}{ludvigak@kth.se}} }
\author{Anna-Karin Tornberg}
\renewcommand\Affilfont{\itshape}
\affil{
  Numerical Analysis, Department of Mathematics,\\
  KTH Royal Institute of Technology, 100 44 Stockholm, Sweden
}
\date{}
\title{\thetitle}
\begin{document}
\maketitle
\fi

\begin{abstract}
  Boundary integral methods are advantageous when simulating viscous
  flow around rigid particles, due to the reduction in number of
  unknowns and straightforward handling of the geometry.
  In this work we present a fast and accurate framework for simulating
  spheroids in periodic Stokes flow, which is based on the completed
  double layer boundary integral formulation. The framework implements
  a new method known as quadrature by expansion (QBX), which uses
  surrogate local expansions of the layer potential to evaluate it to
  very high accuracy both on and off the particle surfaces. 
  This quadrature method is accelerated through a newly developed
  precomputation scheme. The long range interactions are computed
  using the spectral Ewald (SE) fast summation method, which after
  integration with QBX allows the resulting system to be solved in $M
  \log M$ time, where $M$ is the number of particles.
  This framework is suitable for simulations of large particle
  systems, and can be used for studying e.g. porous media models.
\end{abstract}

\ifelsart
\end{frontmatter}
\fi

\section{Introduction}

Fluid flows involving microscopic, rigid particles are common both in
nature and in industrial processes. The macroscopic properties of such
systems are often determined by the particle interactions happening on
the smallest scale of the flow. Examples of such systems can be found
in sedimenting suspensions \cite{Guazzelli2011} and porous media
modeling \cite{DeAnna2013}. To fully understand the interactions
between particles, and the effect those interactions have on the fluid
flow, numerical simulation is a valuable tool. Stokes equations are
often valid in this context, due to the small particle size and low
fluid velocity. For problems governed by Stokes equations it is
possible to use boundary integral methods, where the solution is
represented as a layer potential from the boundaries of the domain (in
this case the particle surfaces).  Since the problem then is
formulated as a boundary integral equation on the union of the
particle surfaces, the dimension of the domain which has to be
discretized is reduced from $\mathbb R^3$ to $\mathbb R^2$,
significantly reducing the number of unknowns. This also removes the
problem of combining a volume grid with moving boundaries, which can
be challenging. However, boundary integral methods come with a set of
challenges of their own, two of which can be considered major.

The first major challenge of boundary integral methods is that the
resulting linear system after e.g. a Nystr\"om discretization is
dense, such that the cost of one left hand side evaluation is
$\ordo(N^2)$, where $N$ is the number of unknowns. This puts a severe
limit on the size of the systems that can be considered, even if a
rapidly converging interative method is used. This can be overcome by
evaluating the layer potential using a fast method, such as the
\ac{FMM} \cite{Greengard1997} or a fast Ewald summation method
\cite{Deserno1998,Lindbo2011c} (for periodic problems). These methods
reduce the cost of one left hand side evaluation to $\ordo(N)$ and
$\ordo(N \log N)$, respectively, thereby making boundary integral
methods suitable for large-scale computations.

The second major challenge of boundary integral methods is that one
needs accurate quadrature of singular and nearly singular integrals
when evaluating the layer potential. Developing such methods which are
accurate, fast and work for arbitrary geometries is a topic of current
research, particularly for the problem of nearly singular quadrature
in three dimensions. In two dimensions there are efficient methods for
nearly singular quadrature which gain a lot of their power from the
complex variable formulation \cite{Helsing2008,Barnett2015}. In three
dimensions the situation appears less resolved, though several
different methods have been successfully used in practical
applications \cite{Bremer2012, Bremer2013, Ying2006, Quaife2014,
  Tlupova2013, Zhao2010}. A common feature of many of these methods is
however that they are highly target specific, meaning that the cost
grows rapidly if there are many nearly singular integrals to be
evaluated.

Quadrature by expansion (QBX) \cite{Klockner2013,Barnett2014} is a
fairly new method for numerical integration of singular and nearly
singular integrals. The method is built around evaluating layer
potentials through local expansions, and comes equipped with a solid
convergence theory \cite{Epstein2013}. It was originally presented for
the Helmholtz kernel in two dimensions, but the principles of the
method generalize directly both to three dimensions and other
kernels. A promising feature of QBX is that it is possible to combine
it with \iac{FMM}, thereby creating a fully $\ordo(N)$ method that is
able to accurately evaluate the potential everywhere.

In this paper we extend the QBX framework to deal with the Stokes
double layer potential in three dimensions, and combine it with the
fast Ewald summation method presented in \cite{AfKlinteberg2014a}. The
result is a robust and scalable framework for computing Stokes flow
around periodic systems of rigid, spheroidal particles. To limit the
scope of the present paper, we restrict our attention to stationary
particles. The method could easily be extended to dynamic problems ---
such as sedimentation --- by coupling it to an ODE solver, as was done
in \cite{AfKlinteberg2014a}.

The structure of this paper is as follows: In section
\ref{sec:formulation} we state the necessary boundary integral
formulation for Stokes flow around rigid particles. In section
\ref{sec:quadrature} we introduce QBX and discuss errors and parameter
selection. In section \ref{sec:accel-qbx-spher} we show how QBX can be
accelerated for our problem, leading to a computationally feasible
method. In section \ref{sec:peri-susp-ewald} we briefly touch on the
subject of Ewald summation, and how to combine it with QBX. Finally,
in section \ref{sec:results} we present selected numerical results,
which we draw both from validation tests and from an example
application. We also include an appendix (\ref{sec:qbx-spheroidal}),
which covers the computational details of the acceleration scheme
presented in section \ref{sec:accel-qbx-spher}.

\section{Formulation}
\label{sec:formulation}

Our application of interest is that of rigid particles in a Newtonian fluid. The particles
are assumed small enough that the Reynolds number is close to zero, such that we can model
the flow as being Stokes flow, governed by the Stokes equations,
\begin{align}
  \begin{split}
    -\nabla P + \mu \Delta \v u + \v f= 0, 
    \\
    \nabla\cdot\v u=0 .
  \end{split}
  \label{eq:stokes_equations}
\end{align}
We here limit ourselves to spheroidal particles with semi-axes $a$ and $c$. The surface of
such a particle, oriented along the $z$ coordinate axis, can in Cartesian coordinates be
described as
\begin{align}
    \frac{x^2+y^2}{a^2} + \frac{z^2}{c^2} = 1.
\end{align}
If parametrized in the spherical coordinates $\theta\in[0,\pi]$ and $\varphi\in[0,2\pi)$,
the same surface is given by
\begin{align}
  \begin{split}
    x &= a \sin(\theta) \cos(\varphi), \\
    y &= a \sin(\theta) \sin(\varphi), \\
    z &= c \cos(\theta).
  \end{split}
  \label{eq:spheroid_param}
\end{align}
The particles can be classified into three distinct subgroups: prolate ($a< c$),
spherical ($a = c$), and oblate ($a> c$).

\subsection{Boundary integral formulation}
\label{sec:bound-integr-form}

We are considering Stokes flow as described by the Stokes equations
\eqref{eq:stokes_equations}, which by virtue of being a set of linear
partial differential equations with constant coefficients have
solutions that can be represented using boundary integrals. For these
representations, the essential Green's functions are the stokeslet
$\stokeslet$, the stresslet $\stresslet$ and the rotlet $\rotlet$%
\footnote{Throughout this paper we use the Einstein convention that an
  index appearing twice in an expression implies a summation over the
  set $\{1,2,3\}$, except for when there is a $\sum$ explicitly
  defining the index.}.
\begin{align}
  \stokeslet_{ij}(\v x,\v y) &= \frac{\delta_{ij}}{r} + \frac{r_ir_j}{r^3}, \\
  \stresslet_{ijk}(\v x,\v y) &= -6 \frac{r_ir_jr_k}{r^5}
  \label{eq:stresslet_def}, \\
  \rotlet_{ij}(\v x,\v y) &= \epsilon_{ijk}\frac{r_k}{r^3} , 
\end{align}
where $\v r = \v x-\v y$ and $r = |\v r|$. 

To compute the flow associated with a rigid particle in a viscous
fluid, we use the completed boundary integral representation of Power
\& Miranda \cite{Power1987}, where the flow in the domain is expressed
as a double layer potential $\v\dblsymb$ plus a completion flow
$\v\ucompsymb$ and a background flow $\ubg$,
\begin{align}
  \v u(\v x) = \dbl{\Gamma,\v q}{\v x} + \ucomp{\v x} + \ubg(\v x).
  \label{eq:u_def}
\end{align}
The double layer potential is the flow associated with a double layer density $\v q$,
defined on the particle surface $\Gamma$. It is computed by integrating $\v q$ over
$\Gamma$ together with the stresslet singularity $\stresslet$ and the outward unit normal
$\v\nhat$,
\begin{align}
  \dbli{\Gamma,\v q}{\v x} = \int_{\Gamma} \stresslet_{ijk}(\v x,\v y) q_j(\v y)
  \nhat_k(\v y) \dSy.
  \label{eq:dbl_def}
\end{align}
The integral in \eqref{eq:dbl_def} is weakly singular for $\v x \in \Gamma$, and is then
taken as the principal value integral, assuming that $\Gamma$ is Lyapunov smooth. Letting
$\v x$ approach the surface from either the interior or exterior domain, there is a jump
in the double layer potential,
\begin{align}
  \lim_{\epsilon\to 0} \dbl{\Gamma,\v q}{\v x \pm \epsilon\v{\nhat}} = 
  \mp 4\pi\v q(\v x) + \dbl{\Gamma,\v q}{\v x}, \quad \v x \in \Gamma.
  \label{eq:dbl_jump}
\end{align}
The completion flow $\v\ucompsymb$ is added to the formulation due to
the double layer potential's inability to represent a net force $\v f$
and torque $\v t$ on the particle \cite{Power1987}. For a rigid
particle, a suitable completion flow is that generated by a stokeslet
singularity $\stokeslet$ of magnitude $\v f/8\pi\mu$ and a rotlet
singularity $\rotlet$ of magnitude $\v t/8\pi\mu$, both located at the
particle center $\v x_c$,
\begin{align}
  \ucompi{\v x} = \frac{1}{8\pi\mu} \left( \stokeslet_{ij}(\v x, \v x_c)f_j + 
    \rotlet_{ij}(\v x, \v x_c)t_j \right).
  \label{eq:ucomp_def}
\end{align}
For a particle undergoing rigid body motion with translational velocity $\v U$ and
rotational velocity $\v\Omega$, we set a no-slip boundary condition,
\begin{align}
  \v u(\v x) = \v U + \v\Omega \times (\v x - \v x_c) .
  \label{eq:noslip}
\end{align}

Letting $\v x$ in our flow field representation \eqref{eq:u_def} go to
the surface of the particle from the exterior, we get a diagonal term
from \eqref{eq:dbl_jump}. We then enforce the no-slip boundary
condition \eqref{eq:noslip} on the surface, and close the system using
\eqref{eq:rbm_V} and \eqref{eq:rbm_Omega}. This results in a \ac{BIE}
of the second kind in the density $\v q$,
\begin{align}
  -4\pi \v q (\v x) + \dbl{\Gamma,\v q}{\v x} + \ucomp{\v x} + \ubg(\v x) = \v U +
  \v\Omega \times (\v x - \v x_c) .
  \label{eq:bie}
\end{align}
The fact that \eqref{eq:bie} is second kind with a compact integral operator is the main
benefit of using a double layer formulation. After discretization using a quadrature rule
and the Nystr\"om method, the corresponding linear system can be said to be well
conditioned in two senses; one is that the condition number is bounded (and usually very
low), the other is that the condition number stays constant under grid refinement (this is
in contrast to single layer formulations, where the condition number increases under
refinement). These properties make the system suitable for iterative solution using
\ac{GMRES}, which typically converges rapidly.

The equation system in \eqref{eq:bie} needs to be closed by specifying either the rigid
body motion $(\v U,\v\Omega)$ or the completion flow $\v\ucompsymb$. This corresponds to
using two different formulations: the resistance problem formulation or the mobility
problem formulation.

\subsubsection{The resistance problem}

When the rigid body motion $(\v U,\v\Omega)$ of the particle is known, the problem is to
find the resulting force $\v f$ and torque $\v t$ exerted by the fluid on the
particle. These can be related to the double layer density $\v q$ through the formulation
by Power \& Miranda \cite{Power1987},
\begin{align}
  \v f(\v q) &= \int_\Gamma \v q(\v y) \dSy, \\
  \v t(\v q) &= \int_\Gamma \v q(\v y) \times (\v y - \v x_c) \dSy,
\end{align}
such that the completion flow is a functional of $\v q$,
\begin{align}
  \ucompi{\v x} = \ucompsymb_i[\Gamma, \v q](\v x) = \frac{1}{8\pi\mu} \left(
    \stokeslet_{ij}(\v x, \v x_c)f_j(\v q) + \rotlet_{ij}(\v x, \v x_c)t_j(\v q) \right).
\end{align}

\subsubsection{The mobility problem}

When the external forcing $(\v f,\v t)$ on the particle is known, the problem is to find
the rigid body motion $(\v U,\v \Omega)$ of the particle. Using the formulation available
in e.g. Pozrikidis \cite{Pozrikidis1992}, the rigid body motion vectors can be computed as
functionals of $\v q$,
\begin{align}
  \v V(\v q) &= -\frac{4\pi}{S_\Gamma} \int_{\Gamma} \v q(\v y) \dSy,   
  \label{eq:rbm_V}
  \\
  \v \Omega(\v q) &= -4\pi
  \sum_{n=1}^3\frac{\gv{\omega}^{(n)}}{A_n}
  \left(
    \v{\omega}^{(n)} \cdot \int_{\Gamma}
    (\v y-\v x_c)\times\v q(\v y)\dSy
  \right),
  \label{eq:rbm_Omega}
\end{align}
where $S_\Gamma$ is the surface area of $\Gamma$ and 
\begin{align}
  A_n = \int_{\Gamma} \left[ \gv{\omega}^{(n)}
    \times (\v y-\v x_c) \right] \cdot \left[ 
    \gv{\omega}^{(n)} \times (\v y-\v x_c) \right] \dSy .
\end{align}
The vectors $\gv{\omega}^{(n)}$, $n=1,2,3$, are independent unit vectors satisfying
\begin{align}
  \frac{1}{\sqrt{A_n A_m}} \int_{\Gamma} 
  \left[ \gv{\omega}^{(m)} \times (\v y-\v x_c) \right] 
  \cdot
  \left[ \gv{\omega}^{(n)} \times (\v y-\v x_c)
  \right] \dSy = \delta_{mn} .
\end{align}

\subsection{Particle systems}
\label{sec:particle-systems}

For a system of $\Npart$ particles, the linearity of the Stokes equations allows us to
express the velocity field as a superposition of the fields from the individual particles,
\begin{align}
  \v u(\v x) = \sum_{\alpha=1}^{\Npart} 
  \left(
     \v\dblsymb\asup(\v x)  + \v\ucompsymb\asup(\v x) 
  \right) + \ubg(v x),
  \label{eq:u_system}
\end{align}
where $\gamalp$ is the surface of particle $\alpha$,
\begin{align}
  \v\dblsymb\asup(\v x) = \dbl{\gamalp,\v q}{\v x},
\end{align}
and $(\v\ucompsymb\asup,\v U\asup,\v\Omega\asup)$ are defined either as in the resistance
formulation or the mobility formulation.

Evaluating \eqref{eq:u_system} and enforcing the no-slip boundary condition
\eqref{eq:noslip} on the surface of each particle gives the boundary integral equation for
a system of particles,
\begin{align}
  \begin{split}    
  -4\pi \v q (\v x) + \sum_{\alpha=1}^{\Npart} \left( \v\dblsymb\asup(\v x) +
    \v\ucompsymb\asup(\v x) \right) + \ubg(\v x)
  = 
  \v U^\beta + \v\Omega^\beta \times (\v x - \v x_c^\beta) ,
  \\
  \v x \in \Gamma_\beta, \quad \beta=1,\dots,M .
  \label{eq:bie_system}
  \end{split}
\end{align}

\subsection{Discrete formulation}

To solve the (single particle) BIE \eqref{eq:bie}, we first need a way of numerically
evaluating integrals over the particle surface, which we denote
\begin{align}
  \opint[f] = \int_\Gamma f(\v y) \dSy .
  \label{eq:def_opint}
\end{align}
This we do using a quadrature rule $\opquad_N$, which defines a set of $N$ nodes $\v x_i$
and weights $w_i$ on $\Gamma$, such that
\begin{align}
  \opint[f] \approx \opquad_{N}[f] = \sum_{i=1}^N f(\v y_i) w_i .
  \label{eq:def_opquad}
\end{align}
The details of the quadrature used will be discussed further in section
\ref{sec:quadrature}. We denote by the superscript $h$ a quantity computed using
$\opquad_N$, e.g.
\begin{align}
  \potential{\numapprox{\dblsymb_i}}{\Gamma,\v q}{\v x} =
  \opquad_N [\stresslet_{ijk}(\v x,\cdot) q_j(\cdot) \nhat_k(\cdot)] .
  \label{eq:dbl_h_def}
\end{align}
Applying the Nystr\"om method, where the integral equation is enforced at the quadrature
nodes, the \ac{BIE} for the single-particle system \eqref{eq:bie} is approximated by the
$3N \times 3N$ linear system in the density values $\v q(\v x_i)$,
\begin{align}
  \begin{split}
    -4\pi \v q (\v x_i) + \dblh{\Gamma,\v q}{\v x_i} + \ucomp{\v x_i} = 
    \v U + \v\Omega \times (\v x_i - \v x_c),
    \\
    \quad i=1,\dots,N,
  \end{split}
  \label{eq:bie_h}
\end{align}
where either $\v\ucompsymb = \numapprox{\v\ucompsymb}$ or $\v U = \numapprox{\v U}$ and
$\v\Omega = \numapprox{\v\Omega}$, depending on if the problem is a resistance problem or
a mobility problem. 

The matrix corresponding to the linear system \eqref{eq:bie_h} is
dense, so the cost of applying a direct solution method to it would be
$\ordo(N^3)$. However, when solving with GMRES the number of
iterations required to reach a certain tolerance is largely
independent of the discretization. The cost is then $\ordo(N^2)$, with
a constant that depends on the problem geometry. This is acceptable
for a single particle, since the solution $\v q$ usually can be
resolved with a limited number of surface nodes. For a problem with
multiple particles, however, the corresponding linear system (derived
from \eqref{eq:bie_system}) is $3\Npart N \times 3\Npart N$ and
dense. The cost then grows quadratically with the number of particles,
which severely limits the size of the problem that can be
considered. The remedy for this is to use a fast method, which
exploits the structure of the problem to evaluate the left hand side
of the problem (also known as the matrix-vector product, or
\emph{matvec}) in less than quadratic time. Common fast methods are
the \ac{FMM} and, for periodic problems, fast Ewald summation. The
respective complexities for these two methods are $\ordo(M)$ and
$\ordo(M \log M)$ as the number of particles is increased. We have
here used fast Ewald summation, which will be further discussed in
section \ref{sec:peri-susp-ewald}.

\section{Quadrature}
\label{sec:quadrature}

To derive the linear system \eqref{eq:bie_h} corresponding to the
\ac{BIE} \eqref{eq:bie}, we needed a quadrature method $\opquad_N$
\eqref{eq:def_opquad} for evaluating the double layer potential
$\v\dblsymb$ \eqref{eq:dbl_def}. We will in this section describe such
a method, which is able to accurately resolve $\v\dblsymb$ also when
the integral is nearly singular or singular.

\subsection{Spheroidal grid}
\label{sec:spheroidal-grid}
When performing quadrature on a spheroidal particle, it is natural to use the parametrized
description \eqref{eq:spheroid_param} in $(\theta,\varphi)$ and to discretize it using an
$\ntheta\times\nphi$ grid in those coordinates, with the grid aspect ratio given by
\begin{align}
  \frac{\nphi}{\ntheta} = \frac{a}{c}.
\end{align}
This gives us a total of $N=\ntheta\nphi$ quadrature points on the surface. We let
$(\theta_i,\lambda_i^\theta)$, $i\in\{1 \dots \ntheta\}$, be the nodes and quadrature weights of
an $\ntheta$-point Gauss-Legendre quadrature rule on the interval $[0, \pi]$, and
$(\varphi_j,\lambda_j^\varphi)$, $j\in\{1 \dots \nphi\}$, be the nodes and quadrature weights of
the trapezoidal rule on the interval $[0, 2\pi)$. This gives us the quadrature rule
\begin{align}
  \opquad_{\ntheta\nphi}[f] = 
  \sum_{i=1}^\ntheta \sum_{j=1}^\nphi
  f(\v y(\theta,\varphi))
  W(\theta_i,\varphi_j)
  \lambda^\theta_i \lambda^\varphi_j
  \label{eq:quad_rule_full}
\end{align}
where $W(\theta,\varphi)$ is the area element. For ease of notation we can simplify this
expression by letting one index cover all quadrature points and introducing $\v
y_i$ %
and $w_i$ %
such that
\begin{align}
  \opquad_{N}[f] = \sum_{i=1}^N f(\v y_i) w_i = \opquad_{\ntheta\nphi}[f].
  \label{eq:quad_rule}
\end{align}
We denote this grid and the accompanying quadrature rule the \emph{spheroidal grid}, and
will refer to the quadrature rule as the direct quadrature.

\begin{figure}
  \centering
  \begin{subfigure}[b]{0.49\textwidth}
    \centering
    \includegraphics[width=0.6\textwidth]{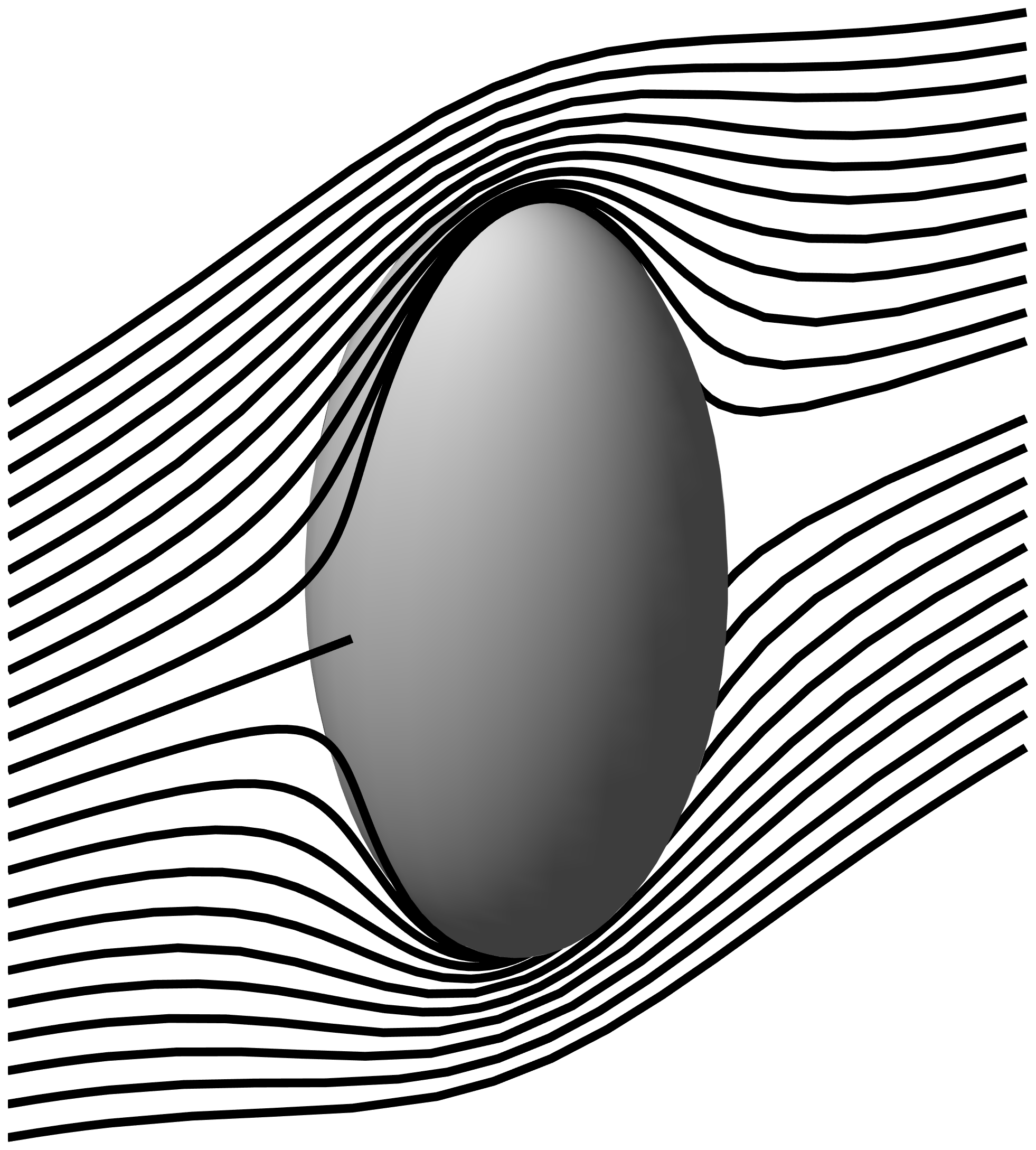}
    \vspace{1.5em}
    \caption{Stokes flow around a single translating spheroid, exact
      solution from \cite{Chwang1975}.}
    \label{fig:stream}
  \end{subfigure}
  \hfill
  \begin{subfigure}[b]{0.5\textwidth}
    \includegraphics[width=\textwidth]{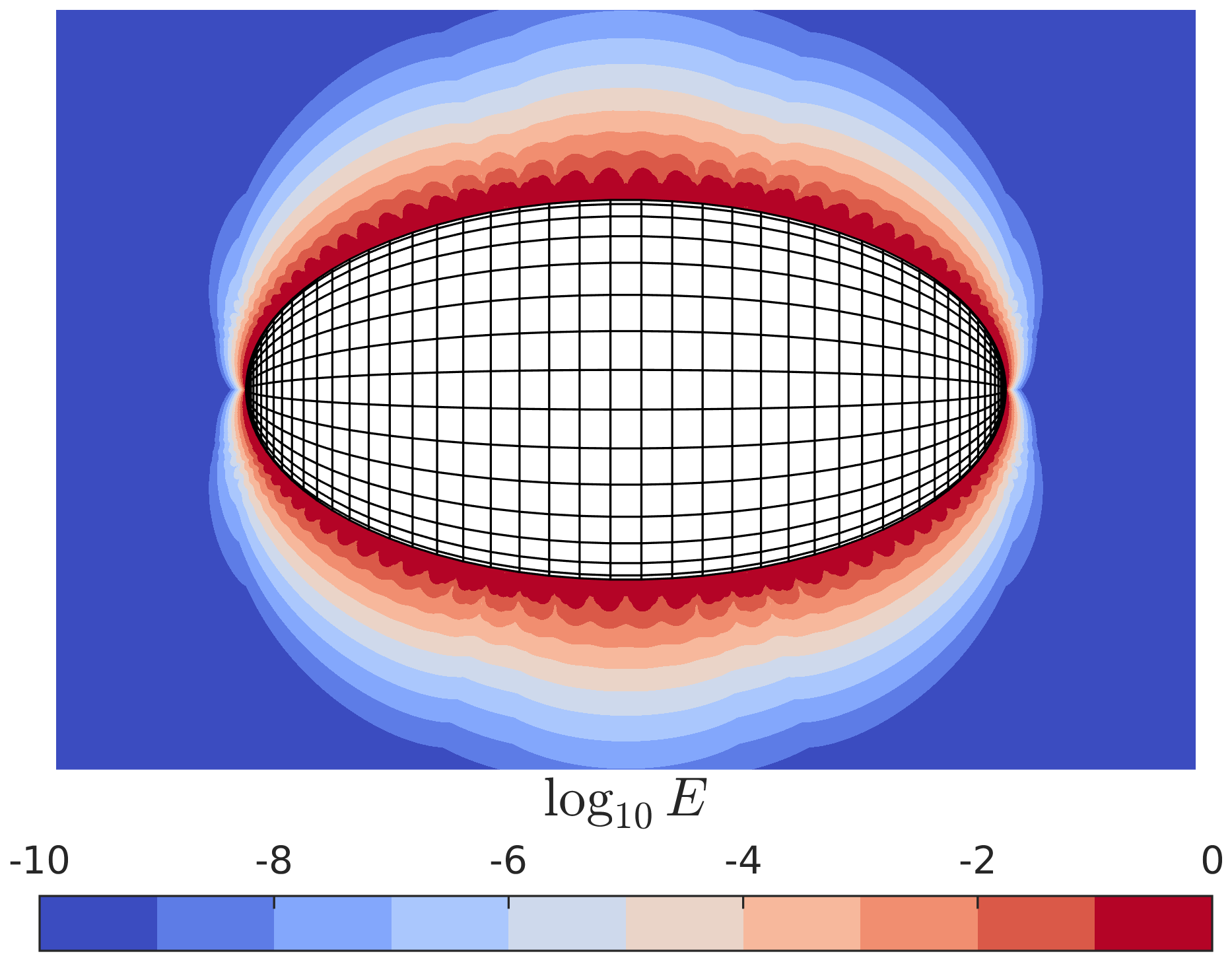}
    \caption{Error when evaluating double layer potential from
      spheroid using direct quadrature.}
    \label{fig:neval_1}
  \end{subfigure}
  \caption{Example of error due to nearly singular integration
    (\subref{fig:neval_1}), when evaluating flow due to a single
    particle (\subref{fig:stream}).}
  \label{fig:quad_err_ex}
\end{figure}

The direct quadrature rule defined in \eqref{eq:quad_rule_full} allows
us to approximate integrals over the surface with spectral accuracy
for \emph{smooth} and \emph{well-resolved} integrands, and works well
for evaluating the double layer potential \eqref{eq:dbl_def} when the
target point $\v x$ is far enough away from the surface
$\Gamma$. However, when $\v x$ is on $\Gamma$ the integral is singular
and can not be computed using a quadrature rule for smooth
functions. There is also a difficulty when $\v x$ lies in the domain,
but is close to $\Gamma$. The singularity in the stresslet then causes
the integral kernel to be sharply peaked, and trying to evaluate it
using a smooth quadrature rule gives an error which grows
exponentially as $\v x$ approaches $\Gamma$ (see Figure
\ref{fig:quad_err_ex}). We refer to this as the integral being
\emph{nearly singular}.

The singular case ($\v x \in \Gamma$) can for the double layer
potential be evaluated using a singularity subtraction method, which
gives third order accuracy at virtually no increased cost
\cite{AfKlinteberg2014a}. For higher order accuracy in the singular
case, and for evaluating the nearly singular case, specialized
quadrature methods are required.

In what follows of this section we will first give a brief
introduction to an existing quadrature method for singular and nearly
singular layer potentials that is based on local expansions. We will
then show how this method can be used to evaluate the Stokes double
layer potential. Lastly, we will show how this method in our case can
be accelerated by geometric considerations.

\subsection{Quadrature by expansion}
\label{sec:qbx}

\begin{figure}[]
  \centering
  \includegraphics{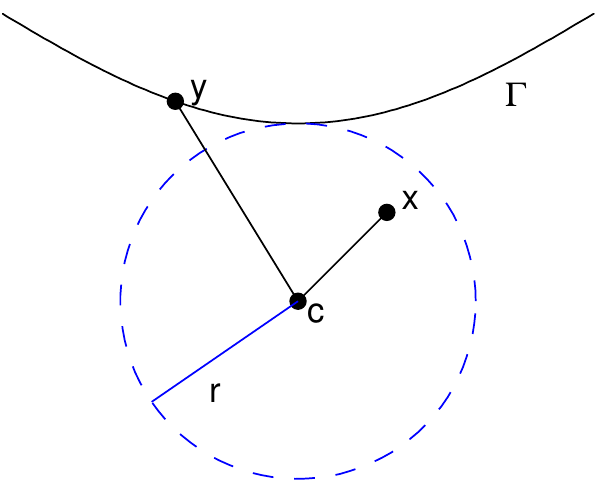}
  \caption{Local expansion and domain of convergence.}
  \label{fig:expansion_geo_curve}
\end{figure}

Quadrature by expansion (QBX) is a new quadrature scheme first
introduced for the Helmholtz equation in two dimensions by Kl\"ockner
et al. \cite{Klockner2013}. The method is based on the observation
that a straightforward quadrature rule applied to a layer potential is
inaccurate when evaluated close to or on the surface, but well
resolved when sufficiently far away, as illustrated in Figure
\ref{fig:quad_err_ex}. Since the layer potential is a smooth function
away from the boundary, we can pick any point $\v c$ in the domain and
create a local expansion of the potential about that point. Picking
that point in the well-resolved part of the domain, we can compute the
coefficients of the expansion using our direct quadrature rule. The
domain of convergence for the local expansion is the ball of radius $r$,
\begin{align}
  r = \min_{\v y \in \Gamma} |\v c - \v y|,
\end{align}
centered on $\v c$, and so the
expansion can be used to accurately evaluate the potential close to
the surface and at the point where the ball touches the surface (see
Figure \ref{fig:expansion_geo_curve} for an example geometry). This is
the essential idea of QBX. For further reading regarding the
convergence and evaluation of these local expansions, we refer to
Epstein et al. \cite{Epstein2013} and Barnett \cite{Barnett2014}.

In three dimensions, a good starting point for the construction of local expansions for
layer potentials is the expansion of the Green's function for the Laplace equation, called
the Laplace expansion,
\begin{align}
  \frac{1}{|\v x-\v y|} = \sum_{l=0}^\infty
  \frac{4\pi}{2l+1} \sum_{m=-l}^l
  r_x^l Y_l^{-m} (\theta_x,\varphi_x)
  \frac{1}{r_y^{l+1}} Y_l^{m}(\theta_y,\varphi_y),
  \label{eq:laplace_expansion}
\end{align}
where $Y_l^m$ is the spherical harmonics function of degree $l$ and order $m$, and
$(r_x,\theta_x,\varphi_x)$ and $(r_y,\theta_y,\varphi_y)$ are the coordinates of $\v x$
and $\v y$ in a spherical coordinate system centered at $\v c$. This expansion gives us
the separation of source ($\v y$) and target ($\v x$) that is necessary for the method. As
a first application, we will consider the double layer potential of the Laplace equation,
which we henceforth will refer to as the dipole potential. For a smooth vector density
$\gv\rho$ defined on the surface $\Gamma$, we define it as
\begin{align}
  \ldbl{\Gamma,\gv\rho}{\v x} := 
  \int_\Gamma \gv\rho\cdot\nabla_{\v y}\frac{1}{|\v x - \v y|} \dSy 
  .
  \label{eq:ldbl_def}
\end{align}
Inserting the expansion \eqref{eq:laplace_expansion} into
\eqref{eq:ldbl_def} and moving the terms related to the target point
$\v x$ out of the integration, we get
\begin{align}
  \ldbl{\Gamma,\gv\rho}{\v x} = 
  \sum_{l=0}^\infty
  \sum_{m=-l}^l
  r_x^l Y_l^{-m} (\theta_x,\varphi_x)
  z_{lm},
  \label{eq:ldbl_expa}
\end{align}
where $z_{lm}$ are the local expansion coefficients of the potential,
\begin{align}
  z_{lm} = \frac{4\pi}{2l+1}
  \int_\Gamma \gv\rho\cdot\nabla 
  \frac{1}{r_y^{l+1}} Y_l^{m}(\theta_y,\varphi_y) \dSy .
  \label{eq:ldbl_exp_coeff}
\end{align}
Once the expansion coefficients $z_{lm}$ are known, they can be used
to evaluate the potential at any target point $\v x$ in the domain of
convergence. At this stage we wish to highlight two things: (i) We
have so far made no approximations, evaluating the potential through
expressions \eqref{eq:ldbl_expa} and \eqref{eq:ldbl_exp_coeff} is
still exact. The sum in $l$ will need to be truncated, and we will
discuss the truncation error that it introduces in section
\ref{sec:error-analysis}. (ii) The integrand in
\eqref{eq:ldbl_exp_coeff} contains no singularity, and can be
evaluated using a smooth quadrature rule. The quadrature errors
related to this smooth quadrature will also be discussed in section
\ref{sec:error-analysis}.

\subsection{QBX for the stresslet}
\label{sec:qbx-stresslet}
To evaluate the Stokes double layer potential \eqref{eq:dbl_def} using
QBX, we need to be able to expand its kernel in a way that separates
source and target, just as we could for the dipole potential. One way
of doing this is to express the double layer potential in terms of the
dipole potential, for which we already have the expansion
\eqref{eq:ldbl_expa}. Here we use the result by Tornberg and Greengard
\cite{Tornberg2008}, which was developed with the purpose of using a
harmonic \ac{FMM} to compute the flow due to a collection of
stokeslets and stresslets. Denoting by $\dipole$ the kernel of the
dipole potential,
\begin{align}
  \dipole_i(\v x, \v y) = \nabla_{\v y}\frac{1}{r} = \frac{r_i}{r^3}, \quad \v
  r = \v x - \v y,
\end{align}
the kernel of the double layer potential can be expressed as
\begin{align}
  \begin{split}
  \stresslet_{ijk}(\v x, \v y)\nhat_k = 
  & \left( 
    (x_j-y_j)\pd{}{x_i} - \delta_{ij}
  \right)
  \dipole_k(\v x, \v y) \nhat _k
  \\
  + & \left(
    (x_k-y_k)\pd{}{x_i} - \delta_{ik}
  \right)
  \dipole_j(\v x, \v y) \nhat_k .    
  \end{split}
\end{align}
Rearranging terms, the integrand of \eqref{eq:dbl_def} can then be written
\begin{align}
  \begin{split}
  \stresslet_{ijk}(\v x,\v y) q_j \nhat_k = 
  &\left(
    x_j\pd{}{x_i} - \delta_{ij}
  \right) \dipole_k(\v x, \v y)(q_j\nhat_k + \nhat_jq_k) 
  \\
  & - \pd{}{x_i} \dipole_j(\v x, \v y)(y_kq_k\nhat_j + y_k\nhat_kq_j).    
  \end{split}
\end{align}
Consequently, the double layer potential can be expressed in terms of dipole potentials as
\begin{align}
  \begin{split}
  \dbli{\Gamma,\v q}{\v x} =
  & \left( x_j\pd{}{x_i}-\delta_{ij} \right) \ldbl{\Gamma,q_j
    \v\nhat+\nhat_j \v q}{\v x}
  \\
  & - \pd{}{x_i}\ldbl{\Gamma,y_kq_k\v\nhat+y_k\nhat_k\v q}{\v x}    
  \end{split}
\label{eq:stresslet_decomp}
\end{align}
With this result we can create an expansion of the double layer
potential by creating four expansions of the dipole potential with
four different densities (including the implicit summation over
$j$). This means that the three vector components of the double layer
potential are computed using four scalar expansions.The expansion
coefficients that need to be computed for a given expansion center $\v
c$ are then
\begin{align}
  \begin{split}
    z_{lm}^j &= 
    \int_\Gamma (q_j \v\nhat+\nhat_j \v q) \cdot\nabla 
    \frac{1}{r_y^{l+1}} Y_l^{m}(\theta_y,\varphi_y) \dSy, 
    \quad j=1,2,3, \\
    z_{lm}^4 &= 
    \int_\Gamma (y_kq_k\v\nhat+y_k\nhat_k\v q) \cdot\nabla 
    \frac{1}{r_y^{l+1}} Y_l^{m}(\theta_y,\varphi_y) \dSy.  
  \end{split}
  \label{eq:exp_coeff}
\end{align}
The spherical harmonics are conjugate symmetric in $m$, $Y_l^{-m}(\theta,\varphi) =
Y_l^m(\theta,\varphi)^*$, and since the density $\v q$ is real the same conjugate symmetry
holds for the expansion coefficients,
\begin{align}
  z_{l,-m}^j = \left( z_{lm}^j \right)^*, \quad j=1,\dots,4 .
\end{align}
Hence, we only need to compute the coefficients for $0\le m \le l$. Truncating the
expansion at some $p = l_{\max}$, the total number of coefficients that we need to compute
for each component $j$ is
\begin{align}
  \Nmom = \frac{p^2 + 3p + 2}{2} .
  \label{eq:def_Nmom}
\end{align}
\subsection{Numerical scheme}

To solve the discrete \ac{BIE} for a single particle \eqref{eq:bie_h},
we need to evaluate the double layer potential at the quadrature
points $\v x_1,\dots,\v x_N$. A local expansion of the potential is
only convergent at a single point on the surface, where the ball of
convergence touches the surface, so we need $N$ expansion centers $\v
c_1,\dots,\v c_N$ located at a distance $r$ normally from the
quadrature points,
\begin{align}
  \v c_i = \v x_i + r \v \nhat_i .
  \label{eq:c_exp}
\end{align}
The same expansions are used to evaluate the potential on the surface
and in the domain. A natural choice of expansion radius is then
\begin{align}
    r= \frac{ d_{\tol} }{2},
\end{align}
where $d_\tol$ has been determined as the minimum
distance from the surface at which the direct quadrature can be used
for a given error tolerance $\tol$. For target points $\v x$ such that
$\min_{\v y\in\Gamma}|\v x-\v y| < d_\tol$, the potential is then
evaluated using the closest expansion center.

\subsection{Error analysis}
\label{sec:error-analysis}

Here we give a sketch of an error analysis for the dipole potential,
with the argument that it carries over to the double layer potential,
as the latter is constructed as a combination of the former. For more
in-depth discussions we refer to \cite{Epstein2013} and
\cite{AfKlinteberg2016quad}.

When evaluating the dipole potential \eqref{eq:ldbl_def} using QBX, we
truncate the expansion \eqref{eq:ldbl_expa} at some order $p$ and use
our direct Gauss-Legendre/trapezoidal quadrature rule
\eqref{eq:quad_rule_full} to approximate the coefficients $z_{lm}$
\eqref{eq:ldbl_exp_coeff}. Denoting the approximate coefficients
$\tilde z_{lm}$, the error is
\begin{align}
  E = 
  \left|
  \ldbl{\Gamma,\gv\rho}{\v x}
  -
  \sum_{l=0}^p \sum_{m=-l}^l
  r_x^l Y_l^{-m} (\theta_x,\varphi_x)
  \tilde z_{lm}
  \right|
  .
  \label{eq:ldbl_err_def}
\end{align}
Adding and subtracting the exact coefficients $z_{lm}$ to the expansion, the error can be
split into two components,
\begin{align}
  E \le E_T  + E_Q,
\end{align}
where $E_T$ is the truncation error due to using a finite number of terms in the
expansion,
\begin{align}
  E_T = \left| \ldbl{\Gamma,\gv\rho}{\v x} - \sum_{l=0}^p \sum_{m=-l}^l r_x^l Y_l^{-m}
    (\theta_x,\varphi_x) z_{lm} \right|,
  \label{eq:def_trunc_error}
\end{align}
and $E_Q$ is the quadrature error due using a discrete quadrature rule,
\begin{align}
  E_Q = \left | \sum_{l=0}^p \sum_{m=-l}^l r_x^l Y_l^{-m} (\theta_x,\varphi_x) (\tilde
    z_{lm} - z_{lm}) \right| .
  \label{eq:def_quad_error}
\end{align}
From the results in \cite{Epstein2013} we have that the truncation error can be expected
to behave as
\begin{align}
  E_T = \ordo(r^{p+1}) .
\end{align}
In \cite{AfKlinteberg2016quad} an estimate was derived for the QBX quadrature error, when
considering the Laplace single layer potential on a spheroid. Simplifying those results,
we can expect that the quadrature error should behave approximately as
\begin{align}
  E_Q = \ordo\left( \left( \beta \frac{r}{h} \right)^p e^{-\alpha r/h} \right),
\end{align}
for some $\alpha,\beta > 0$ and $h$ a measure of the grid size; for the spheroidal grid we
define
\begin{align}
  h := 2\pi a/n_\varphi .
  \label{eq:h_def}
\end{align}
Putting the two expression together, we see that the total error has one component that
that decays with $p$ and one component that grows with $p$,
\begin{align}
  E = \ordo(r^{p+1}) + \ordo\left( \left( \beta \frac{r}{h} \right)^p e^{-\alpha r/h}
  \right).
\end{align}
This means that if we fix $r$ and $h$ for a given problem, then there exists an optimal
value for $p$ that minimizes the error. If we instead fix the relation $r/h$ at some
constant value, e.g. $r=3h$, then the quadrature error stays constant ($E_Q=C_Q$) such
that for a given $p$
\begin{align}
  E = \ordo(h^{p+1}) + C_Q .
  \label{eq:error_order}
\end{align}
The scheme can then be viewed as having order of accuracy $(p+1)$ under grid refinement,
up to the point where the truncation and quadrature errors are equal. After that we have
to refine the grid (while maintaining a constant $r$) to get any additional improvements.

In practice we want to be able to get arbitrarily high accuracy out of
QBX, without having to introduce more degrees of freedom. To overcome
the problem of the quadrature error growing with $p$, we therefore
\emph{upsample} the surface density to a factor $\kappa$ finer grid
before computing the coefficients. The fine grid is of the same type
as the original, but with $\kappa\ntheta \times \kappa\nphi$ points,
and the density on the fine grid is computed using trigonometric
interpolation and barycentric Lagrange interpolation
\cite{Berrut2004}. This allows us to write the total error as
\begin{align}
  E = \ordo(r^{p+1}) + \ordo\left( \left( \beta \frac{\kappa r}{h} \right)^p
    e^{-\alpha\kappa r/h} \right).
\end{align}
The truncation and quadrature errors can then be controlled
separately, since the exponential term will rapidly reduce the
quadrature error when $\kappa$ is increased.

\begin{figure}
  \centering
  \begin{subfigure}[b]{0.48\textwidth}
    \includegraphics[width=\textwidth]{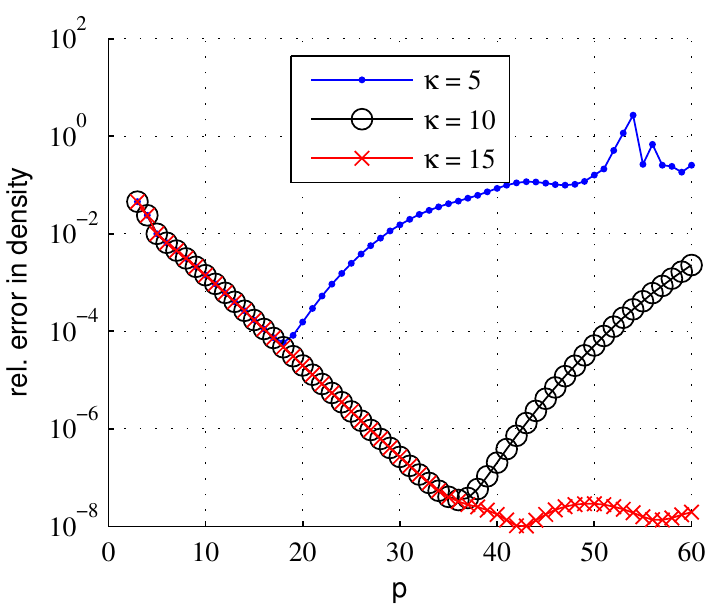}
    \caption{$L^2$ error in the density when solving the mobility
      problem on an oblate spheroid ($a/c=2/1$) with a $40\times 20$
      grid using $r/h=3/2$ and a range of $p$ and $\kappa$.}
    \label{fig:p_conv}    
  \end{subfigure}
  \hfill
  \begin{subfigure}[b]{0.48\textwidth}
    \includegraphics[width=\textwidth]{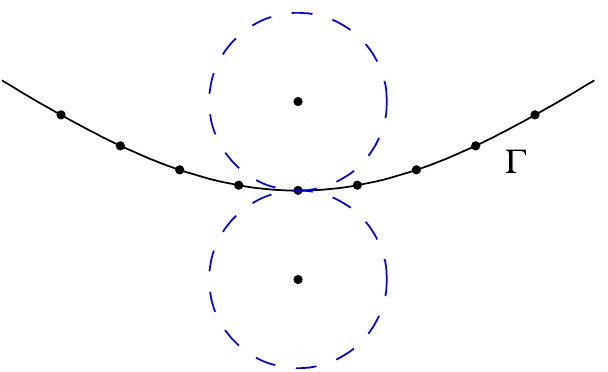}
    \vspace{1em}
    \caption{Expansion centers are used on both sides of $\Gamma$, to
      accurately capture the principal value integral.\\~}
    \label{fig:expansion_geo_twoside}    
  \end{subfigure}
  \caption{ }
\end{figure}

To illustrate the effects of the parameters $p$ and $\kappa$, we solve the mobility
problem for a single oblate spheroid ($a/c=2$) with a $20 \times 40$ grid using $r/h=3/2$
and a range of $p$ and $\kappa$. Measuring the $L^2$ norm of the error in the density
compared to a reference solution, we can see in Figure \ref{fig:p_conv} that the error
decreases exponentially in $p$ until it hits a floor given by the resolution of the
underlying grid, but only for $\kappa$ large enough. If $\kappa$ is taken too small, the
error increases almost exponentially after reaching a minimum. Finding the optimal
$\kappa$ for a given $p$ is not trivial, and picking $\kappa$ larger than necessary is
costly, since it adds a factor $\kappa^2$ to the computational complexity.

\subsection{Two-sided expansions}

When evaluating the singular double layer potential on a surface, we want the principal
value integral. However, what we get when evaluating it using QBX is the one-sided limit
from the side of the expansion center, which we denote $\dblsymb^+$ or $\dblsymb^-$
depending on side,
\begin{align}
  \v\dblsymb^\pm[\Gamma,\v q](\v x) = 
  \lim_{\epsilon\to 0} \dbl{\Gamma,\v q}{\vx \pm \epsilon \v\nhat}.
\end{align}
The relation between the principal value and the one-sided limit is
known through the jump relation \eqref{eq:dbl_jump}, which can be used
together with the QBX result to recover the principal value. However,
and this was originally noted in \cite{Klockner2013}, it turns out
that the spectrum of the linear system better approximates that of the
original operator if we compute the principal value through the
relation
\begin{align}
  \v\dblsymb = \frac{\v\dblsymb^+ + \v\dblsymb^-}{2} .
\end{align}
To do this in practice, we need to put one expansion center on each
side of the target point on the surface and take the average potential
from the two to get the principal value, as illustrated in Figure
\ref{fig:expansion_geo_twoside}. This seemingly doubles the amount of
work needed, but through the acceleration scheme of section
\ref{sec:accel-qbx-spher} that extra cost is hidden from the
computations.

\subsection{Parameter selection}
\label{sec:parameter-selection}

\begin{figure}
  \centering
  \begin{subfigure}[t]{0.49\textwidth}
    \includegraphics[width=\textwidth]{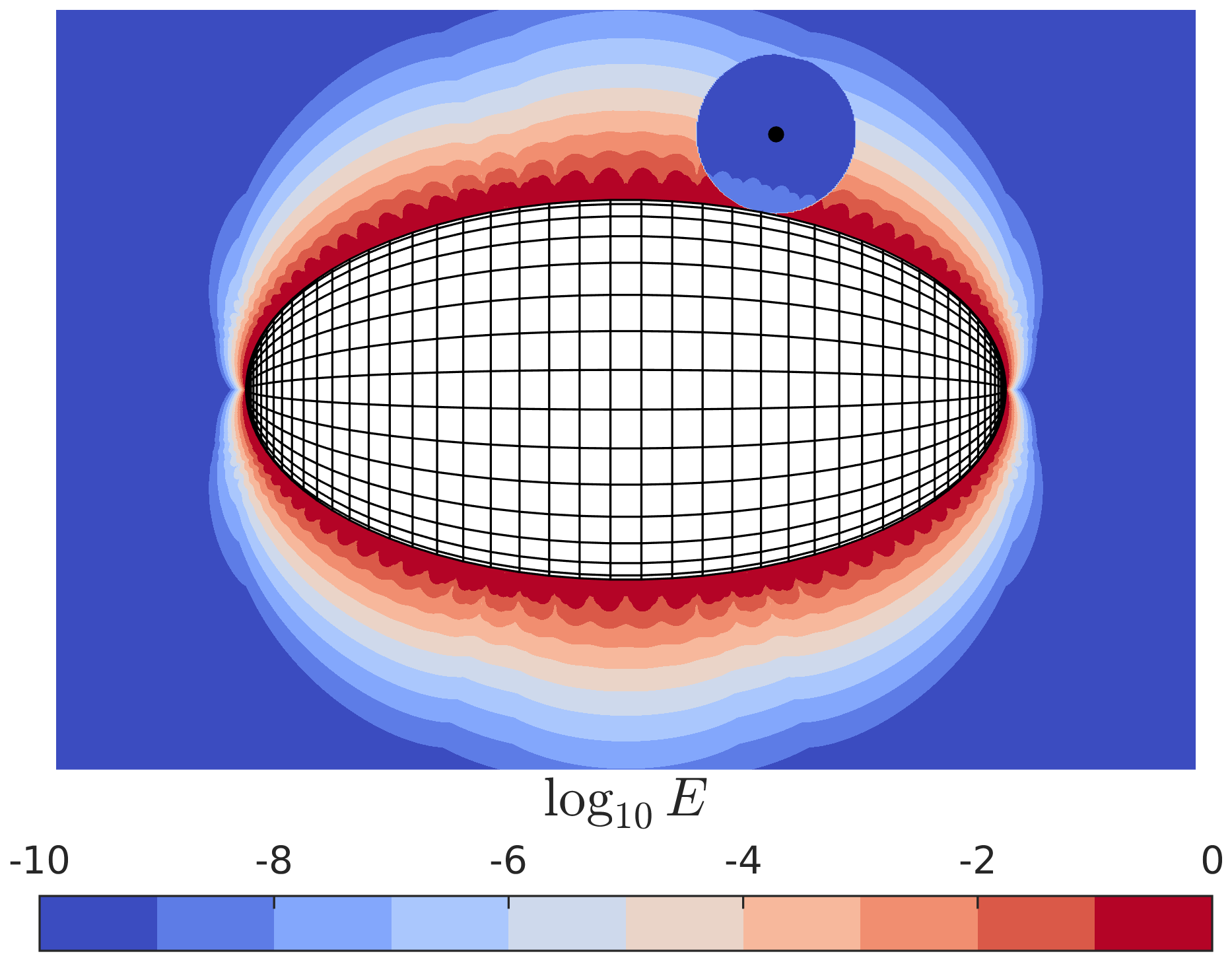}
    \caption{Error when using direct quadrature, punctured in disc
      around an expansion center.}
    \label{fig:neval_punc}
  \end{subfigure}
  \hfill
  \begin{subfigure}[t]{0.49\textwidth}
    \includegraphics[width=\textwidth]{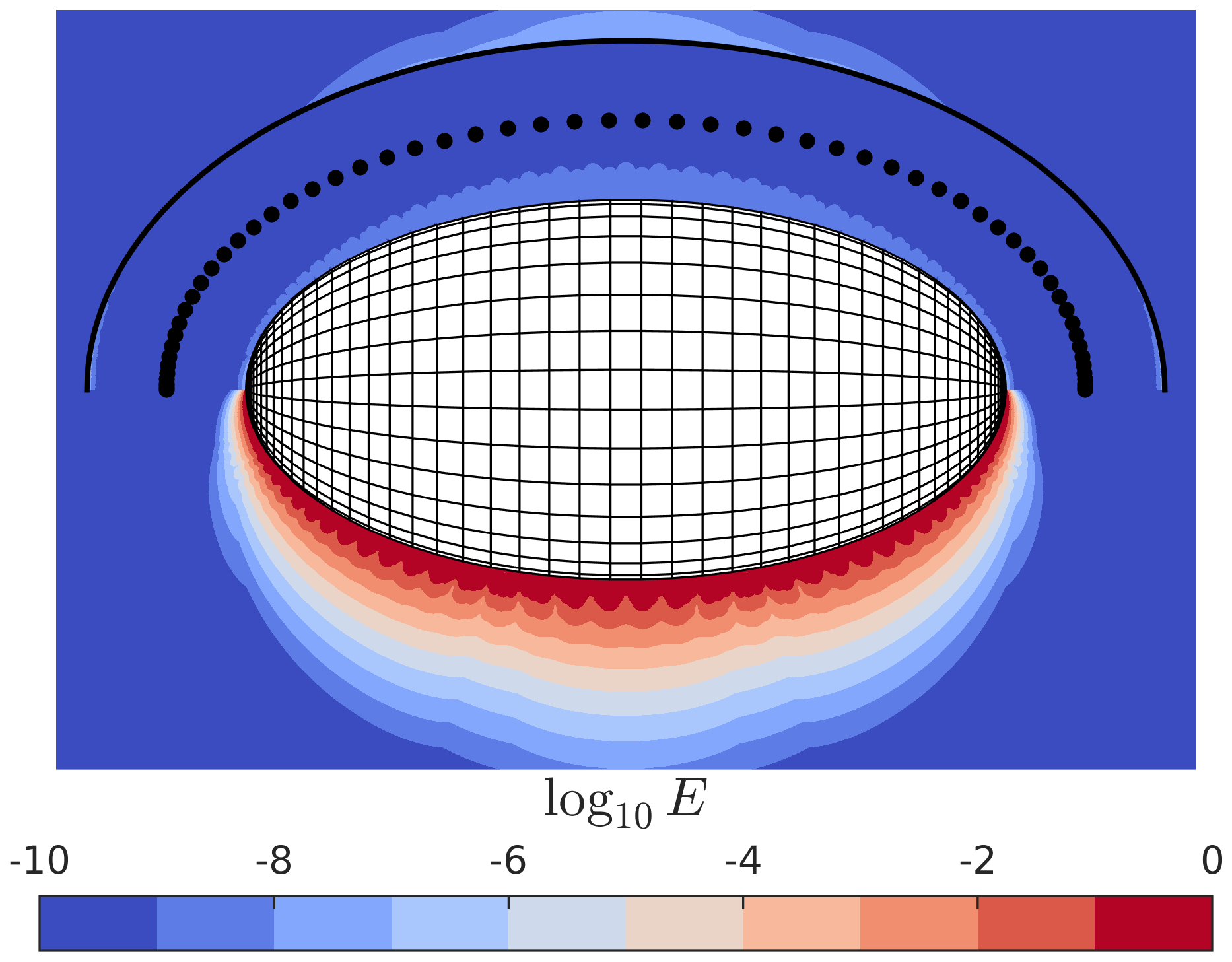}
    \caption{Error when QBX is used in top half plane. Solid black
      line marks threshold $d_\tol$.}
    \label{fig:neval_qbx}
  \end{subfigure}
  \caption{Relative error in solution around prolate spheroid,
    $(a,c)=(1,2)$, compared to exact solution from
    \cite{Chwang1975}. Grid is $60\times 30$, QBX parameters are
    $p=40$, $k=10$, $r/h=2$. Note the error pattern inside the QBX
    area, which is due to an error in the solution, not the
    quadrature.}
  \label{fig:neval}
\end{figure}

For a given geometry, grid resolution $h$ and tolerance $\tol$, an empirical process for
selecting the QBX parameters $p$, $\kappa$ and $r$ is as follows:
\begin{enumerate}
\item Set $r = d_\tol/2$, such that the error from the regular
  quadrature will be less than epsilon when the evaluation point is at
  least a distance $d_\tol$ from $\Gamma$.
\item Set $p$ large enough such that the truncation error $E_T$ from
  the QBX expansion is smaller than $\tol$ for points closer than
  $d_\tol$ from $\Gamma$.
\item Set $\kappa$ large enough to match $p$, such that the quadrature
  error $E_Q$ is kept below $\tol$.
\end{enumerate}
Figure \ref{fig:neval} shows an example where the QBX parameters have
been tuned for high accuracy.

\section{Acceleration of QBX on spheroidal grid}
\label{sec:accel-qbx-spher}

Evaluating the potential from one particle using the QBX scheme
described in section \ref{sec:quadrature} involves a lot of work. For
a spheroidal grid with $N$ points, a given upsampling factor $\kappa$
and expansion order $p$, we need to interpolate the density to a
factor $\kappa^2$ larger grid and then compute $\ordo(p^2)$ expansion
coefficients at every expansion center. This puts the cost for
computing the potential from a particle on itself, the
''self-interaction'', at $\ordo(\kappa^2p^2N^2)$. One natural way of
speeding up the process is to incorporate the QBX scheme into
\iac{FMM}, which hierarchically computes potentials using expansions
of the same kind as those used in QBX (the method was in fact
conceived with this in mind). This makes it possible to compute
everything on the fly at a computational cost that is linear in $N$
(but still includes $\kappa^2p^2$). The modifications are however
nontrivial, and the only unified FMM/QBX method published to date (by
Rachh et al. \cite{Rachh2016}) is for two-dimensional problems.

Another way of achieving speedup is to precompute as much as possible
of the work associated with QBX, and then use precomputed values
whenever possible. This comes at a potentially large storage cost, and
still has a computational complexity that is quadratic in $N$, albeit
with a significantly smaller leading constant. In this work we have
used precomputations in combination with a fast Ewald summation
method, which computes the long-range interactions between $M$
particles in $\ordo(M \log M)$ time (more on this in section
\ref{sec:peri-susp-ewald}). This is feasible partly because the number
$N$ of grid points required on a spheroid is limited, such that the
overall computational cost scales like $\ordo(M \log M)$, and partly
because the spheroidal geometry of the particles in our problem
significantly reduces the storage required for the precomputed
values. We will in this section outline the principles and main
results of our precomputation strategy; the details of the
implementation are available for reference in appendix
\ref{sec:qbx-spheroidal}.

\subsection{Self-interaction}

When solving the boundary integral equation \eqref{eq:bie} on a single
particle using the spheroidal grid of section
\ref{sec:spheroidal-grid}, we typically use the Nystr\"om method to
create a linear system which we solve iteratively. We then have $N$
grid points denoted $\v x_1, \dots, \v x_N$ on the surface, and for
every iteration we need to compute the self-interaction at those
points. With $\v Q \in \mathbb{R}^{3N}$ containing $\{\v q (\v
x_i)\}_{i=1}^N$, we can let $R_i \in \mathbb{R}^{3 \times 3N}$ be the
map $\v Q \to \v u(\v x_i)$ using QBX, such that
\begin{align}
  \v u(\v x_i) = \mself_i \v Q .
  \label{eq:self_interaction}
\end{align}
The matrix $\mself_i$ then represents the resulting action after
upsampling the density to a fine grid, using the fine grid to compute
the expansion coefficients at the expansion center $\v c_i$, and
finally evaluating that expansion at $\v x_i$. Precomputing all the
matrices $\mself_1,\dots,\mself_N$ requires $3 N\times 3N$ memory
storage, but once computed they allow the computation of $\v Q \to \v
u(\v x_i)$ at all points $\v x_i$ at a the cost of a matrix-vector
product, i.e. $\mathcal O(N^2)$. This means that the $\mathcal
O(p^2\kappa^2)$ cost related to the QBX parameters --- and ultimately
to the accuracy --- only appears in the precomputation step. Once the
$R_i$ matrices are computed the self-interaction evaluation has a
fixed cost, independent of which QBX parameters were used.

The memory cost of precomputing all the self interaction matrices can be reduced by taking
advantage of the fact that the spheroidal grid is rotationally symmetric about the polar
axis. This symmetry means that the matrices $\mself_a$ and $\mself_b$ related to two
points $\v x_a$ and $\v x_b$ on the same latitude ($\theta_a = \theta_b$) can be related
as
\begin{align}
  \mself_a = T_z \mself_b (T_z^{-1} \otimes P),
\end{align}
where $T_z \in \mathbb{R}^{3\times 3}$ represents a rotation about the
polar axis of the particle, and $P \in \mathbb{R}^{N \times N}$ is a
permutation matrix, i.e. a row permutation of the identity matrix that
represents a permutation of the point ordering. The symbol $\otimes$
denotes the Kronecker product. As a consequence, it is enough to
precompute matrices related to the $\ntheta$ points on the first
longitude, and then use rotations and permutations to compute the
potential at all grid points with those matrices.

\begin{figure}
  \centering
  \begin{subfigure}[b]{0.45\textwidth}
    \centering
    \includegraphics{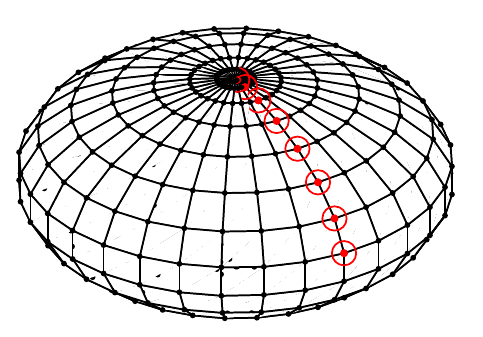}
    \vspace{1em}
    \caption{Spheroidal grid with the $\ntheta/2$ first points circled.}
    \label{fig:spheroid_comp_points}
  \end{subfigure}
  \hfill
  \begin{subfigure}[b]{0.45\textwidth}
    \centering
    \includegraphics[width=0.55\textwidth]{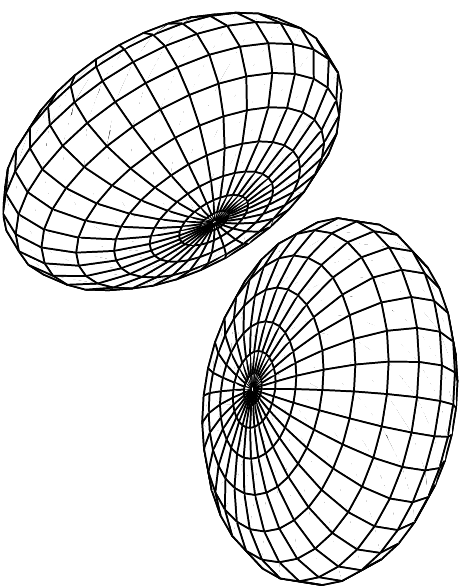}
    \caption{Example of configuration where near-singular quadrature is necessary.}
    \label{fig:two_close_spheroids}
  \end{subfigure}
  \caption{}
  \label{fig:spheroid_acc}
\end{figure}

The savings in memory and precomputation time can be further increased by taking advantage
of the fact the spheroidal grid has a mirror symmetry about the equator. For two points
$\v x_b$ and $\v x_c$ on the same longitude ($\varphi_b=\varphi_c$) but different sides of
the equator, such that $\theta_b = \pi-\theta_c$, we can relate the matrices $\mself_b$
and $\mself_c$ as
\begin{align}
  \mself_b = T_{x} R_c (T_{x}^{-1} \otimes F),
\end{align}
where $T_x \in \mathbb{R}^{3\times 3}$ represents a rotation that mirrors the points $\v
x_b$ and $\v x_c$ into each other, and $F \in \mathbb{R}^{N \times N}$ is a permutation
matrix. Taking advantage of this mirror symmetry allows a further reduction of the
required number of precomputed matrices to $\ntheta/2$, corresponding to the points marked
in Figure \ref{fig:spheroid_comp_points}.

\subsection{Nearly singular evaluation}

When particles are close to each other, such as in Figure \ref{fig:two_close_spheroids},
we need to use QBX to evaluate their interaction accurately, since at least some of the
quadrature points on one particle will be too close to the other particle for the direct
quadrature to be accurate.  To be able to evaluate the potential at arbitrary points close
to a particle using QBX, we need the expansion coefficients for the potential at all the
expansion centers in the domain outside the surface.

We let $\v z^j(\v c_i) \in \mathbb{C}^\Nmom$, $j \in \{1 \dots 4\}$,
denote the four vectors of $\Nmom$ \eqref{eq:def_Nmom} coefficients
each at a an expansion center $\v c_i$, as defined in
\eqref{eq:exp_coeff}. Analogously to \eqref{eq:self_interaction}, the
process of upsampling the density and computing the expansion
coefficients at all the expansion coefficients can then by represented
by a set of matrices $\mmom^j_i \in \mathbb{C}^{\Nmom \times 3N}$,
such that
\begin{align}
  \begin{split}
    \begin{array}{ll} 
          \v z^j(\v c_i) = \mmom^j_i\v Q, \quad &i=1,\dots,N, \\
          &j=1,\dots,4 .
    \end{array}    
  \end{split}
\end{align}
Using the same arguments of symmetry as above, we can relate these matrices by rotational
symmetry as
\begin{align}
  M_a^j = \sum_{i=1}^4 A_{ij} E M_b^i (T_z^{-1} \otimes P),
\end{align}
and by mirror symmetry as
\begin{align}
  M_b^j = \sum_{i=1}^4 B_{ij} \Theta^j M_c^j (T_x^{-1} \otimes F) ,
\end{align}
where $\{\Theta,E\} \in \mathbb{R}^{\Nmom \times \Nmom}$ are diagonal, and
$\{A,B\}\in\mathbb R^{4\times 4}$ orthogonal. Just as for the matrices $\mself_i$, it is
therefore sufficient to precompute the matrices $\mmom_i^j$ for the first $\ntheta/2$
expansions points, since the remaining matrices can be reconstructed by symmetry
operations.

The precomputation of the matrices $\mmom_i^j$ hides the cost
$\kappa^2$ related to the upsampling from the nearly singular
evaluation, though the $\ordo(p^2)$ expansion order cost still affects
the evaluation and storage cost. This is because it reflects the
number of expansion coefficients used, and we need all the
coefficients to able to evaluate the expansion at arbitrary target
points. Only when precomputing for a specific target point, as in the
self-interaction, can we get rid of the $p^2$ factor.

\subsection{Storage complexity}
Precomputing the density-to-coefficient matrices $\mmom_i^j$ for a single expansion center
$\v c_i$ requires storage of $12{\Nmom}N$ complex numbers, while the density-to-potential
matrix $\mself_i$ related to a target point $\vx_i$ only requires storage of $9N$ real
numbers. Considering a spheroidal body, where $N=\nphi \ntheta$, and taking both
symmetries into account, the storage required for computing all necessary $\mmom_i^j$ is
then $6{\Nmom}\ntheta^2 \nphi$ complex numbers, while the storage required for all
necessary $\mself_i$ is $\frac{9}{2} \ntheta^2 \nphi$ real numbers.

As an example, using a set of $50 \times 50$ quadrature points and $p=15$ requires
approximately 1.6 GB for the $\mmom_i^j$ and 4.5 MB for the $\mself_i$, using double
precision. The number of expansion coefficients scales as ${\Nmom} = \mathcal O(p^2)$, so
assuming $\nphi \propto \ntheta = n$ allows us to write the scaling of these
expressions as $\mathcal O(p^2n^3)$ and $\mathcal O(n^3)$. Alternatively, introducing a
grid size $h\propto 1/n$ gives $\mathcal O(p^2h^{-3})$ and $\mathcal O(h^{-3})$.

\subsection{Summary of precomputation scheme}
The precomputation scheme outlined above minimizes the cost for QBX when computing it
directly, as opposed to accelerating it by \iac{FMM}. For a given setup of particle shape,
quadrature grid and QBX parameters, we can precompute the matrices $\mself_i$ and
$\mmom_i^j$, associated with self-interaction and nearly singular evaluation, and store
them to disk. The symmetries of the spheroidal grid allow us to do this at a manageable
cost in terms of memory. Once we have these precomputed values, we can load and use them
for simulations with a large number of particles that all have identical shape and
grid. The single set of precomputed values can then be used for all particles in the
simulation, since they will differ from each other only by a rigid body rotation.

For the self-interaction, the costs of upsampling and expansion order
are hidden from the evaluation step by the precomputations. We can
therefore pick $\kappa $ and $p$ large enough to make sure that that
the accuracy in the evlauation of the double layer potential is only
limited by how well the density $\v q$ is resolved on the underlying
quadrature grid. For the nearly singular evaluation the expansion
order $p$ still affects the computational cost, so there is still a
cost/accuracy tradeoff there. However, the upsampling factor $\kappa$
is hidden by the precomputations, so we can always make sure that we
pick it large enough to resolve integrands associated with large $p$.

\section{Periodic systems and Ewald summation}
\label{sec:peri-susp-ewald}

To be able to model a large particle system without boundary effects, we use the common
method of applying periodic boundary conditions. We then consider a system of $\Npart$
particles contained in a box of size $L_1 \times L_2 \times L_3$, which is periodically
replicated in all three spatial directions. The velocity field is then periodic,
\begin{align}
  \v u^P(\v x) = \v u^P(\v x+\v\tau) ,
  \label{eq:u_per}
\end{align}
where the superscript $P$ denotes the periodicity and $\v\tau$ is a periodic shift,
\begin{align}
  \v\tau = (p_1 L_1, p_2 L_2, p_3 L_3), \quad \v p \in \mathbb{Z}^3 .
  \label{eq:tau_def}
\end{align}
The periodic velocity field is computed as a superposition of the field from all periodic
replications of all particles,
\begin{align}
  \v u^P(\v x) = \sum_{\alpha=1}^{\Npart} \left[ \dblsymb^P[\gamalp,\v q](\v
    x) + \ucompsymb^{P,\alpha}(\v x) \right],
  \label{eq:u_system_per}
\end{align}
where $\dblsymb^P$ is the periodic double
layer potential,
\begin{align}
  \dblsymb^P_i[\Gamma,\v q](\v x) = 
  \sum_{\v\tau} 
  \int_{\Gamma} \stresslet_{ijk}(\v x,\v y+\v\tau) q_j(\v y) \nhat_k(\v y) \dSy,
  \label{eq:dbl_def_per}
\end{align}
and $\ucompsymb^{P,\alpha}$ is the periodic completion flow,
\begin{align}
  \ucompsymb_i^{P,\alpha}(\v x) = \sum_{\v\tau} \ucompsymb_i^\alpha (\v x + \v\tau).
  \label{eq:ucomp_def_per}
\end{align}
Computing $\v u^P$ by direct summation over all periodic images is not
feasible due to the sums over $\v\tau$ in \eqref{eq:dbl_def_per} and
\eqref{eq:ucomp_def_per} only being conditionally convergent. Here we
compute them using a fast Ewald summation method, which we have also
used in previous work \cite{AfKlinteberg2014a}. The stresslet
$\stresslet$ in \eqref{eq:dbl_def_per} is then decomposed into two
parts,
\begin{align}
  \stresslet(\v x, \v y) = \stresslet\rspace(\v x, \v y, \xi) + \stresslet\kspace(\v x, \v
  y, \xi),
  \label{eq:stresslet_decomp_ewald}
\end{align}
where $\stresslet\rspace$ is sharply peaked and short-ranged, while $\stresslet\kspace$ is
smooth and long-ranged. The Ewald parameter $\xi$ determines how short-ranged and smooth
the two respective parts are. Away from the source $\v y$ the short-ranged part
$\stresslet\rspace$ decays exponentially,
\begin{align}
  \stresslet\rspace(\v x, \v y, \xi) \sim e^{-\xi^2|\v x-\v y|^2}.
\end{align}
This makes it possible to introduce a truncation radius $r_c$ and only consider
near-neighbor interactions within that radius ($|\v x - \v y| \le r_c$). The long-ranged
part $\stresslet\kspace$ converges rapidly when computed in Fourier space, where its
Fourier coefficients $\widehat{\stresslet}\kspace(\v k, \xi)$ decay exponentially due to
its smoothness,
\begin{align}
  \widehat\stresslet\kspace(\v k, \xi) \sim e^{-k^2/4\xi^2} .
\end{align}
The Fourier expansion is truncated at a maximum wave number $K = \max_{\v k} \|\v k
\|_\infty$ and computed using the \ac{SE} method \cite{Lindbo2010,Lindbo2011c}, which uses
\iac{FFT} to compute the long-range interactions to spectral accuracy. In
\cite{AfKlinteberg2014a} we describe in detail how this is done for the stresslet, and
provide accurate truncation error estimates which can be used for selecting the Ewald
parameter $\xi$ and the cutoffs $r_c$ and $K$. The equivalent development for the
stokeslet and rotlet singularities are available in \cite{Lindbo2010} and
\cite{AfKlinteberg2016rot}. Computing $\stresslet\kspace$ using the \ac{SE} method and
only computing $\stresslet\rspace$ for near neighbors gives a method with $\ordo(N \log
N)$ complexity, where $N$ is the total number of degrees of freedom when the number of
particles in the system is increased, under the assumption that the average number of near
neighbors of each discrete point is kept constant \cite{AfKlinteberg2014a}.

\subsection{QBX and Ewald summation}
\label{sec:qbx-ewald-summation}

In section \ref{sec:quadrature} we developed a framework for computing the double layer
potential \eqref{eq:dbl_def} to high accuracy for both the singular and nearly singular
cases, based on an expansion for the stresslet $\stresslet = r_i r_j r_k / r^5$. When
computing the periodic double layer potential \eqref{eq:dbl_def_per} using Ewald
summation, we have split the stresslet into a short-ranged part $\stresslet\rspace$ and a
long-ranged part $\stresslet\kspace$ \eqref{eq:stresslet_decomp_ewald}. The explicit form of the
short-ranged part is \cite{AfKlinteberg2014a}
\begin{align}
  \begin{split}
    \stresslet_{jlm}\rspace(\v x, \v y, \xi) =& -\frac{2}{r^4}\left( \frac{3}{r}\erfc(\xi r) 
      + \frac{2\xi}{\sqrt\pi}(3+2\xi^2 r^2-4\xi^4 r^4)e^{-\xi^2r^2} \right) 
    r_j r_l r_m \\
    &+ \frac{8\xi^3}{\sqrt\pi}(2-\xi^2 r^2)e^{-\xi^2r^2} (\delta_{jl}
    r_m+\delta_{lm}r_j+\delta_{mj} r_l),
  \end{split}
  \label{eq:ewald_real}
\end{align}
where $\v r = \v x - \v y$ and $r = |\v r|$. All of the singular behavior of the stresslet
is contained in $\stresslet\rspace$, since a series expansion of $\stresslet\rspace$ in
the limit of $r \to 0$ shows that
\begin{align}
  \lim_{\v x \to \v y} \left[ \stresslet(\v x, \v y) - \stresslet\rspace(\v x, \v y,
    \xi)\right] = 0.
\end{align}
This means that we need a way of using QBX for evaluating the short-ranged potential,
which we define as
\begin{align}
  \dblsymb\rspace_i[\Gamma,\v q](\v x) = 
  \int_{\Gamma} \stresslet_{ijk}\rspace(\v x,\v y,\xi) q_j(\v y) \nhat_k(\v y) \dSy .
  \label{eq:dbl_rspace}
\end{align}
Deriving a local expansion for $\stresslet\rspace$, as defined in
\eqref{eq:ewald_real}, is a daunting proposition. Instead, we rewrite
the short-ranged potential as
\begin{align}
  \dblsymb\rspace_i[\Gamma,\v q](\v x)
  = &\int_{\Gamma} \stresslet_{ijk}(\v x,\v y) q_j(\v y) \nhat_k(\v y) \dSy \\
  & -\int_{\Gamma} %
  \underbrace{ \left[ %
      \stresslet_{ijk}- \stresslet_{ijk}\rspace %
    \right] }_{\stresslet_{ijk}\kspace} %
  (\v x,\v y,\xi) q_j(\v y) \nhat_k(\v y) \dSy .
  \label{eq:dbl_rspace_decomp}
\end{align}
The first term is now the ordinary double layer potential \eqref{eq:dbl_def}, for which we
can use our existing QBX framework. The second term is the smooth and long-ranged part of
the potential, which we can evaluate using our direct quadrature rule. In our
implementation we use \eqref{eq:dbl_rspace_decomp} together with QBX if $\v x$ is close to
or on $\Gamma$, and \eqref{eq:dbl_rspace} truncated at $r_c$ otherwise.

\subsection{Computational complexity}

When now consider a system of $\Npart$ particles with $N$ degrees of freedom each, which
we solve using our accelerated QBX scheme combined with \ac{SE}. The computational
complexity of solving this system then has three distinct components: (i) the interactions
of particles with themselves (the self-interactions), (ii) the nearly singular
interactions between particles, and (iii) the well-separated interactions. The
self-interactions are completely precomputed, and can therefore be evaluated at a
$\ordo(\Npart N^2)$ cost. The nearly singular interactions use the precomputed
$M_i^j$-matrices, and contribute with a $\ordo(p^2\Npart N^2)$ cost, $p$ being the QBX
expansion order and the leading constant depending strongly on the interparticle
spacing. The well-separated interactions are handled by \ac{SE}, at a cost which scales as
$\ordo(\Npart N \log \Npart N)$ as we increase the size of the particle system. This puts
the total computational complexity at
\begin{align}
  \ordo(\Npart N^2) + \ordo(p^2\Npart N^2) + \ordo(\Npart N \log \Npart N) .
\end{align}
For a given particle system with fixed $N$, $p$ and average
interparticle spacing, the cost then grows as $\ordo(\Npart \log
\Npart)$ as the system size is increased by increasing the number of
particles $M$.

\section{Numerical results}
\label{sec:results}

\begin{figure}[htbp]
  \centering
  \begin{subfigure}[b]{0.49\textwidth}
    \centering
    \includegraphics[width=.95\textwidth]{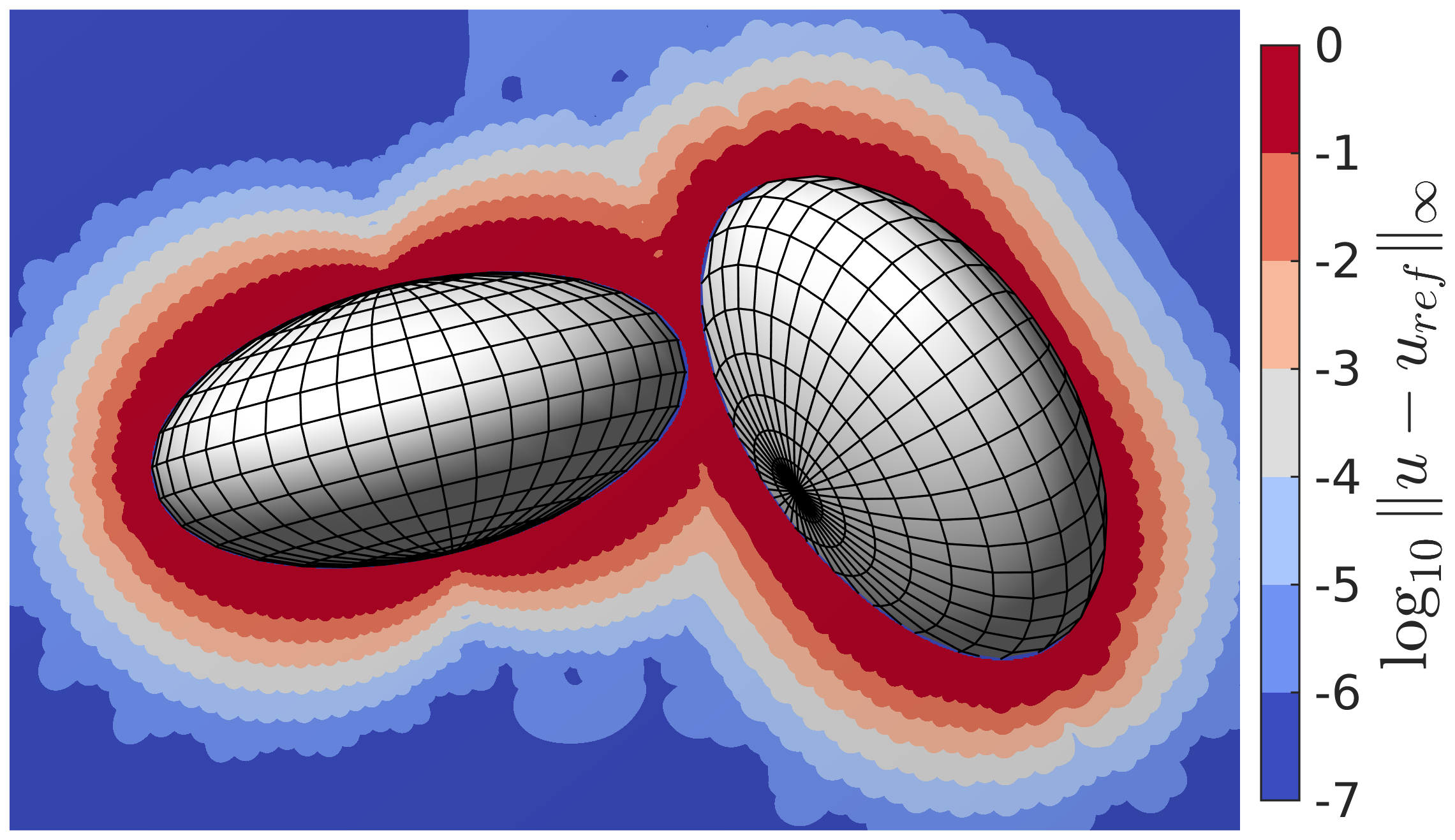}
    \caption{Error using direct quadrature.}
  \end{subfigure}
  \begin{subfigure}[b]{0.49\textwidth}
    \centering
    \includegraphics[width=.95\textwidth]{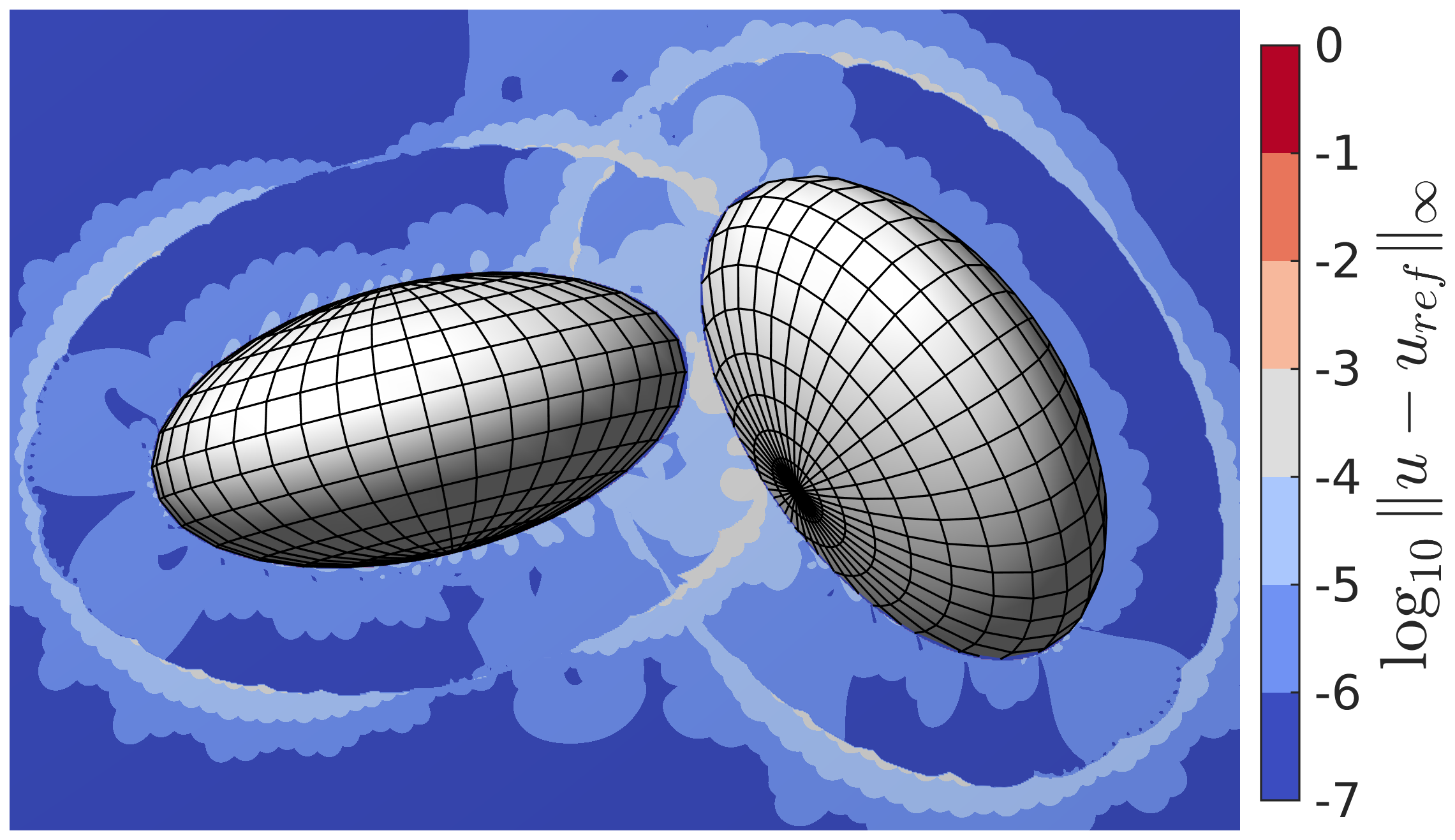}
    \caption{Error using QBX.}
  \end{subfigure}
  \caption{Example of quadrature error (vs. reference solution) when
    evaluating the flow around two close particles, with and without
    QBX. The particles are oblates, $a/c=2$, with a $20 \times 40$
    grid and QBX parameters $r/h=1.5$, $\kappa=10$ and $p=30$.}
  \label{fig:twoparticle_err}
\end{figure}

We can now summarize our numerical scheme as follows:
\begin{itemize}
\item \textbf{Spheroidal grid:} For a spheroidal body with axes $a$
  and $c$, we introduce a spheroidal grid (section
  \ref{sec:spheroidal-grid}) of dimensions $N=\ntheta \times
  \nphi$. For a system with $M$ particles, an identical grid is used
  for each particle.
\item \textbf{QBX:} We place expansion centers $\v c_i$ \eqref{eq:c_exp} at a normal
  distance $r$ from the grid points $\v x_i$, and use an order $p$ local expansion
  (section \ref{sec:qbx}) to evaluate the double layer potential \eqref{eq:dbl_def} in the
  vicinity of each expansion center. The expansion is computed using a factor $\kappa$
  upsampled grid. The QBX computations are accelerated by the precomputation scheme of
  section \ref{sec:accel-qbx-spher}.
\item \textbf{Ewald:} For periodic problems, we compute the velocity field by combining
  QBX with a fast and spectrally accurate Ewald summation method (section
  \ref{sec:peri-susp-ewald}).
\item \textbf{Solution:} We solve the discrete system \eqref{eq:bie_h}
  with $3MN$ unknowns iteratively using GMRES, with the left hand side
  computed using QBX and Ewald.
\item \textbf{Preconditioning:} For systems with many particles, we construct a simple
  block-diagonal preconditioner by computing the explicit inverse of the corresponding
  single-particle system and using it to create the diagonal blocks (through
  rotation). This requires $(3N)^2$ memory storage, but reduces the number of
  iterations by as much as a factor three for some problems.
\end{itemize}
This scheme allows us to efficiently compute the Stokes flow around
systems of particles where the particles are very close to each other,
and also to accurately evaluate the flow everywhere in the domain. An
example is shown in Figure \ref{fig:twoparticle_err}, where one can
clearly see the benefit using QBX when evaluating the potential close
to the particles.

\subsection{Validation}

To validate that our method is correct, we have compared it against
several analytical, semi-analytical and numerical results. For a single
particle in free space, analytic results are available in Chwang \& Wu
\cite{Chwang1975,Chwang1975a} and Jeffery \cite{Jeffery1922}. Figure
\ref{fig:neval} shows an example of a comparison against
\cite{Chwang1975}. For systems of three spheres in free space,
semi-analytic results for several configurations are available in
Wilson \cite{Wilson2013}. For a periodic array of spheres, numeric
results are available in Zick \& Homsy \cite{Zick1982} and Sangani \&
Acrivos \cite{Sangani1982}. The described method can solve all these
problems to arbitrary accuracy, provided that the grid is fine enough
to represent the solution, and that the QBX parameter selection of
section \ref{sec:parameter-selection} is followed. We will below
discuss results for a three-sphere system and for a periodic array of
spheres.

Our implementation of the method has also been used in a separate work
by Bagge \cite{Bagge2015}, where it was thoroughly validated against
the solutions by Chwang \cite{Chwang1975a} and Jeffery
\cite{Jeffery1922}.

\subsubsection{Triangle of spheres}
\label{sec:triangle-spheres}

\begin{figure}[htbp]
  \centering
  \includegraphics[width=0.49\textwidth]{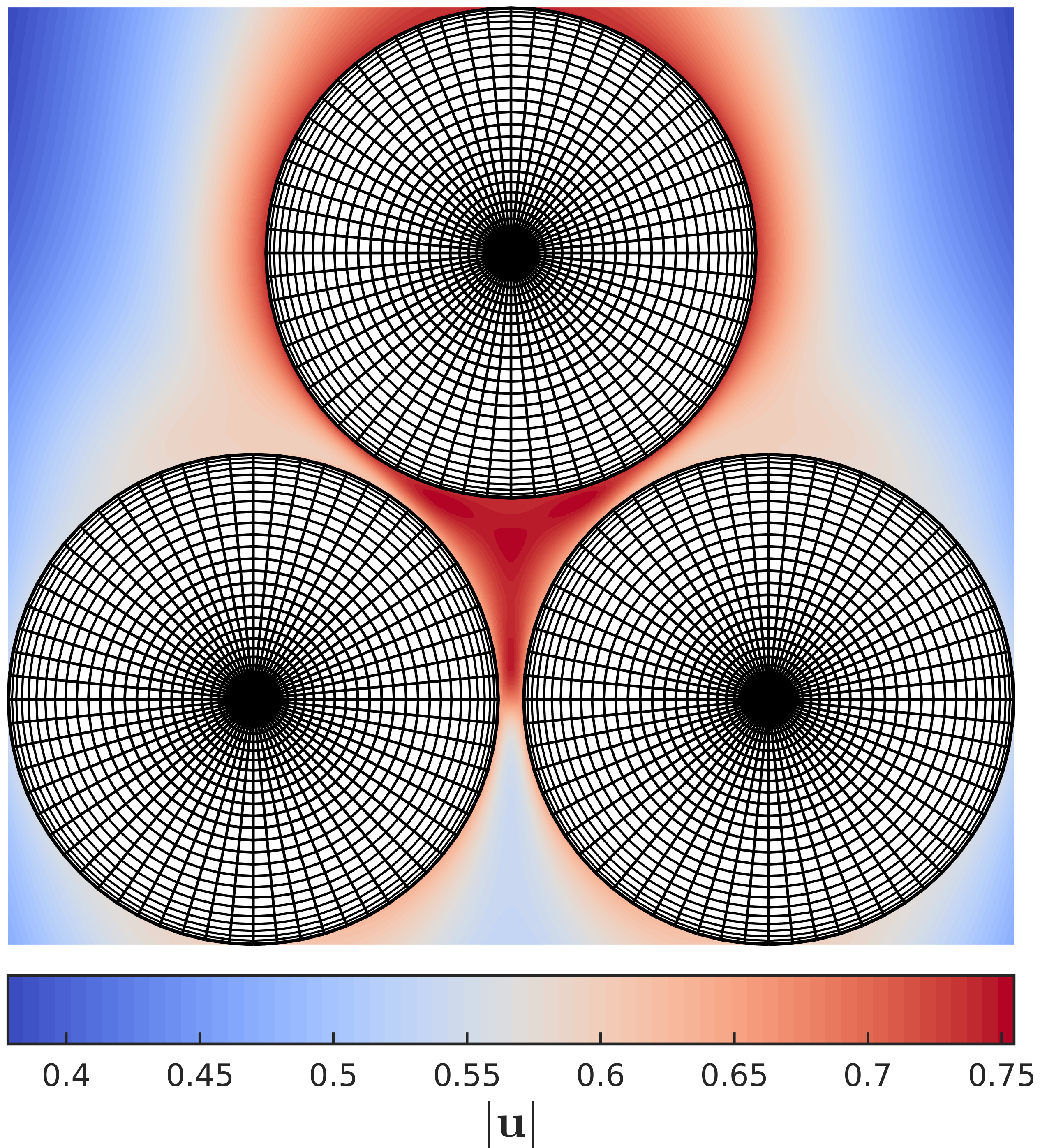}
  \hfill
  \includegraphics[width=0.49\textwidth]{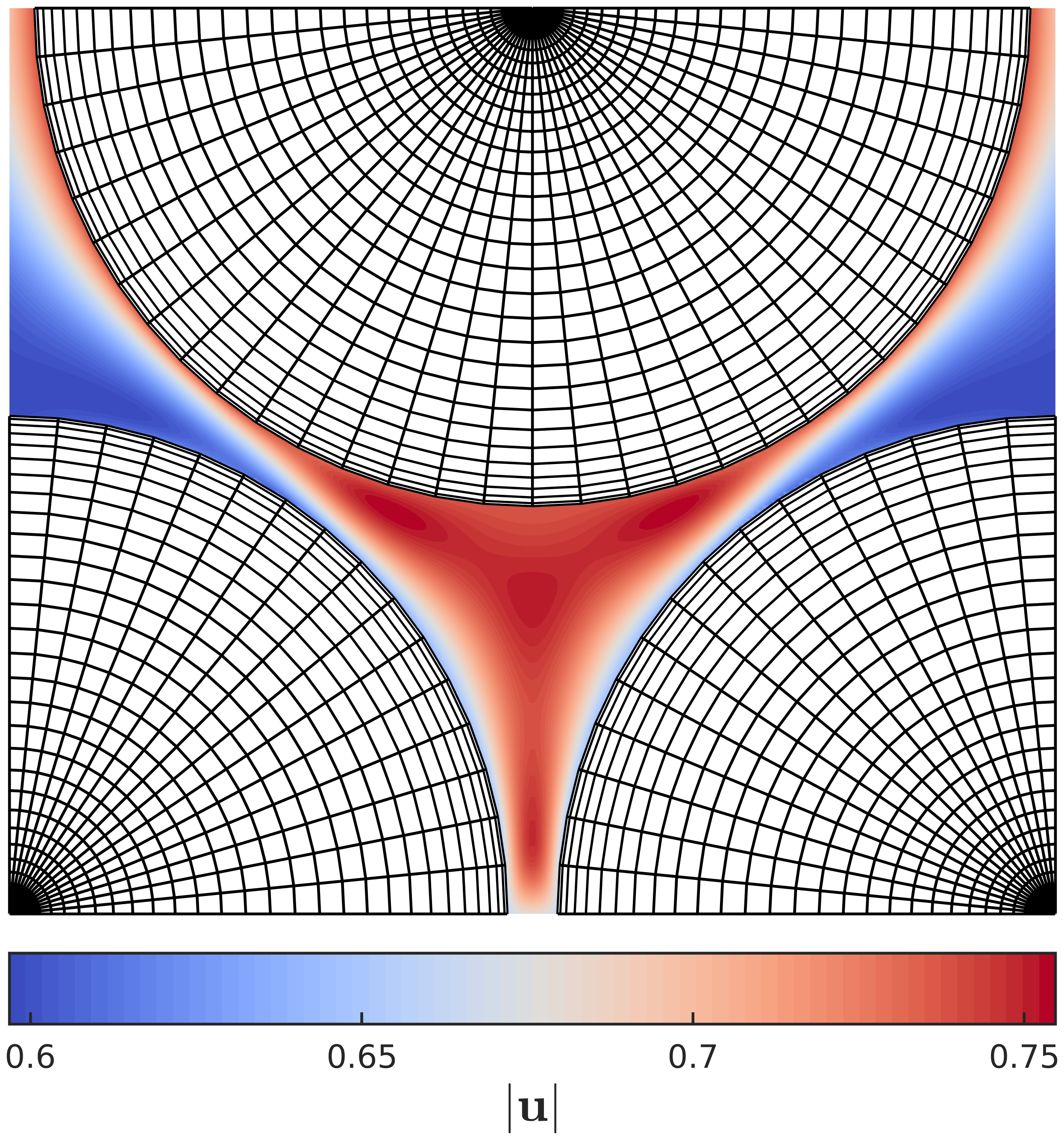}
  \caption{Flow velocity around triangle of spheres at configuration
    $s=2.1$. Sphere 1 is the top sphere. A feature of boundary
    integral methods is that once the solution $\v q$ is known, it is
    possible to zoom in and study the details of the flow field
    everywhere, as is done here in the right picture.}
  \label{fig:triangle}
\end{figure}

To test the ability of our method to resolve particle interactions in
the lubrication limit, we set up the following test case from Wilson
\cite{Wilson2013}: Three spheres of radius $a$ are arranged in an
equilateral triangle with sides $sa$, $s>2$. A force of magnitude
$6\mu\pi a$ is applied to sphere 1 towards the centroid of triangle,
pushing the sphere into the other two spheres. An example of this
geometry for $s=2.1$ is shown in Figure \ref{fig:triangle}. The
following quantities are measured:
\begin{itemize}
\item Velocity $U_1$ of sphere 1, in direction of force.
\item Velocity $U_2$ of spheres 2 and 3, in direction of force.
\item Velocity $U_3$ of spheres 2 and 3, in direction perpendicular to force.
\item Angular velocity $a\Omega$ of spheres 2 and 3.
\end{itemize}

\begin{table}[htbp]
  \centering
  \begin{tabular*}{\textwidth}{@{\extracolsep{\fill}}l llll llll}
    \hline \\
      & \multicolumn{4}{l}{Wilson} & \multicolumn{4}{l}{Current method} \\
    \cline{2-5} \cline{6-9}
    s
    & $U_1$ & $U_2$ & $U_3$ & $\Omega$
    & $U_1$ & $U_2$ & $U_3$ & $\Omega$ \\
    \hline
    2.01& 0.65528& 0.63461& 0.00498& 0.037336 & 0.655\underline{01} & 0.634\underline{24} & 0.00\underline{201} & 0.03\underline{6978}\\
    2.10& 0.73857& 0.59718& 0.03517& 0.052035 & 0.73857 &   0.59718 &   0.03517 &   0.05203\underline{6}                            \\
    2.50& 0.87765& 0.49545& 0.07393& 0.045466 & 0.87765 &   0.49545  &    0.07393 &   0.045466\\
    3.00& 0.93905& 0.41694& 0.07824& 0.035022 & 0.93905 &   0.41694  &    0.07824 &   0.035022\\\\
    \hline
  \end{tabular*}
  \caption{
    Translational and angular velocities for triangle of spheres, computed using $64 \times 64$ grid, $r/h=3$, $\kappa=5$ and $p=20$. Deviations from reference results by Wilson \cite{Wilson2013} are underlined.
  }
  \label{tab:wilson}
\end{table}

As can be seen in Table \ref{tab:wilson}, we are able to reproduce the
results by Wilson to high accuracy for $s \ge 2.10$ using a $64 \times
64$ grid. For $s=2.01$ the results are inaccurate even though the
quadrature is accurately resolved everywhere, and a finer grid is
required if the results are to be improved. This is because the double
layer density $\v q$ gets sharply peaked as the particles get close to
each other, with a peak that has spatial scale $\mathcal O(d^{1/2})$
when the separation distance is $d$ \cite{Barnett2015,Sangani1994a},
and we can not hope to get good results if our grid is too coarse to
resolve the density. The problem is therefore one of resolution, not
quadrature. To solve the problem efficiently as $s \to 2$, an adaptive
refinement strategy for the surface grid would be required.

\subsubsection{Periodic array of spheres}

\begin{figure}[htbp]
  \centering
  \begin{subfigure}[b]{0.4\textwidth}
    \centering
    \includegraphics[width=\textwidth]{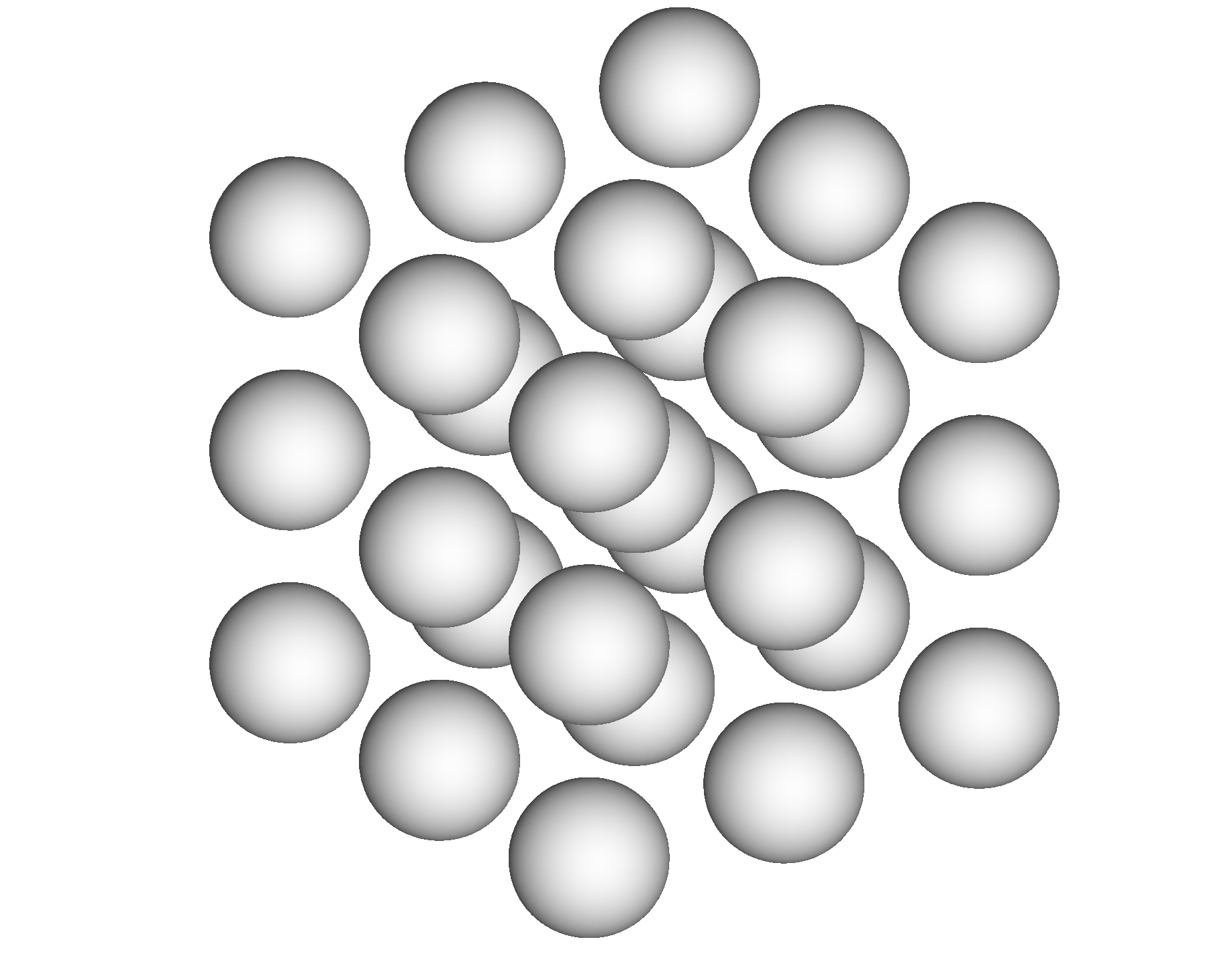}
    \caption{Simple cubic array.}
    \label{fig:sc_lattice}
  \end{subfigure}
  \hspace{0.05\textwidth}
  \begin{subfigure}[b]{0.4\textwidth}
    \centering
    \begin{tabular}{c|c|c}
      $\chi$ & $K$ & $K_\text{ref}$ \\
      \hline
      0.3~ &  ~1.699884 & ~1.7000  \\
      0.4~ &  ~2.151801 & ~2.1518  \\
      0.5~ &  ~2.842022 & ~2.8420  \\
      0.6~ &  ~3.973781 & ~3.9738  \\
      0.7~ &  ~6.004034 & ~6.004~  \\
      0.8~ &  10.054098 & 10.05~~  \\
      0.9~ &  19.161078 & 19.16~~  \\
      0.95 &  27.918287 & 27.9~~~  \\
      1.00 &  42.102343 & 42.1~~~
    \end{tabular}
    \caption{Dimensionless drag.}
    \label{tab:sc_table}
  \end{subfigure}
  \caption{Dimensionless drag coefficient $K$ for a simple cubic array of spheres,
    compared to reference results from \cite{Sangani1982}. The dimensionless parameter
    $\chi=(\rho / \rho_\text{max})^{1/3}$ measures how close the configuration is to the
    touching configuration ($\chi=1$).}
  \label{fig:periodic_array}
\end{figure}

As a validation test for our periodic computations, we consider a
periodic array of spheres of radius $a$ arranged in a simple cubic
lattice. This is modeled using a single sphere in a cubic box with
sides $L$. We subject the sphere to an external force, and solve the
mobility problem to get its velocity for various values of the
concentration $\rho=4\pi a^3 / 3L^3$. We then compute the
dimensionless drag coefficient $K = F / 6 \pi a U$ and compare it to
the results available from Sangani \& Acrivos
\cite{Sangani1982}. Using a $30 \times 30$ grid and tuning all
parameters (QBX+Ewald) for a tolerance of $\tol = 10^{-7}$, we get the
results shown in Figure \ref{fig:periodic_array}. Our results are
identical to the reference results in all the digits presented in
\cite{Sangani1982}, all the way up to the limit $\rho =
\rho_\text{max}$, where the spheres touch ($L=2a$). Our results for
$K$ also agree in all six decimals with the series expansion in
$\chi=(\rho / \rho_\text{max})^{1/3}$ published in
\cite[eq. 60]{Sangani1982}, up to and including $\chi=0.5$, after
which the series expansion loses precision.

\subsection{Application: Porous media flow}

As a demonstration of a possible application for our method, we
consider a periodic cell of dimensions $40 \times 15 \times 15$
containing 131 randomly positioned, oblate particles with size
$(a,c)=(2,1)$, for a volume particle concentration of about 24\%
system can be used as a model of a three-dimensional porous medium
with non-spherical grains. To this end, we set a background flow
$\ubg=(1,0,0)$ and solve the resistance problem for fixed particle
positions. We use a $30 \times 60$ grid and tune our parameters for a
tolerance $\tol = 10^{-4}$. The resulting system has 707,400 unknowns
and is solved in about 2 hours on a regular workstation with a 3.4 GHz
quad-core Intel Haswell CPU.

\begin{figure}[htbp]
  \centering \subcaptionbox{\label{fig:porous_example}}{%
    \includegraphics[width=.62\textwidth]{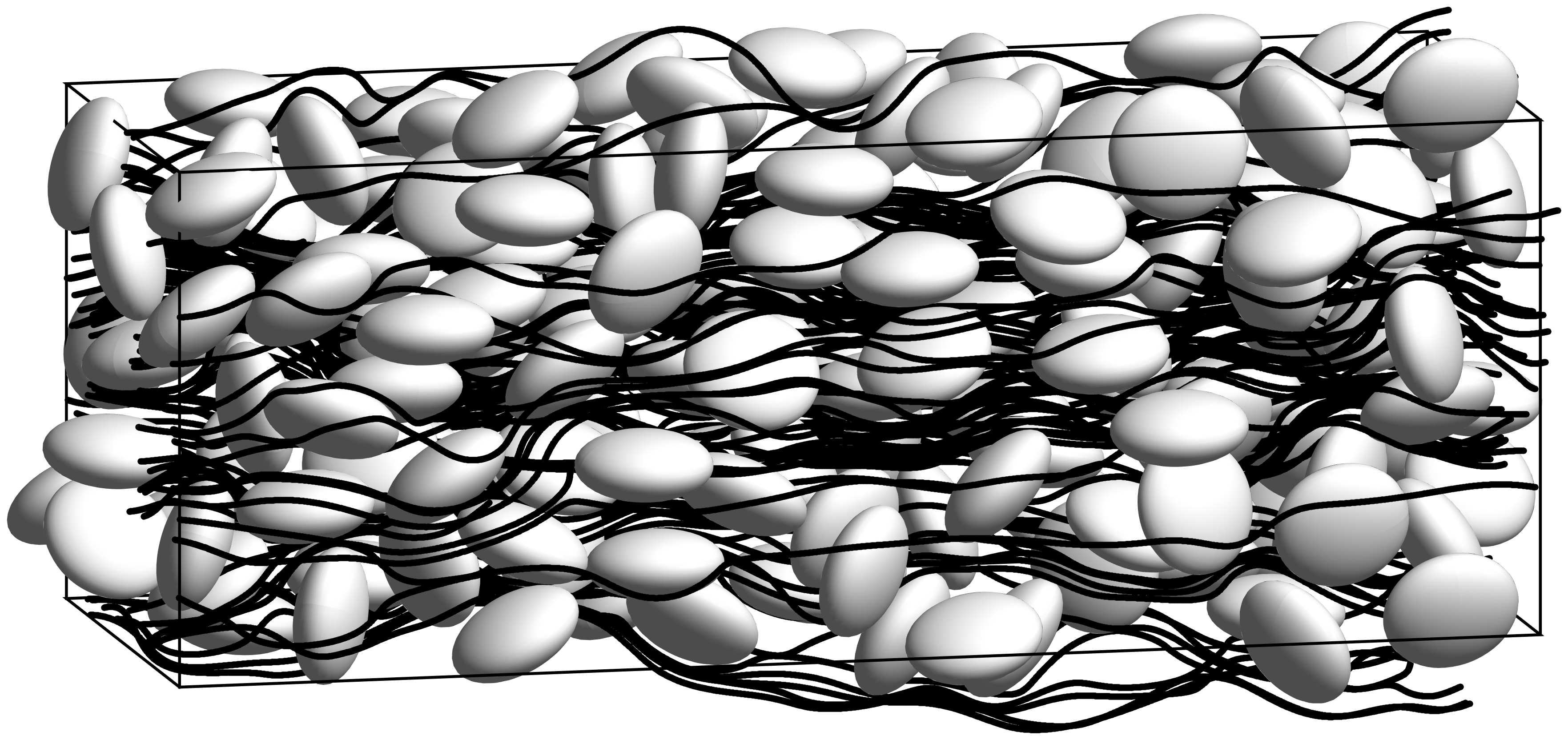}\vspace{.07\textwidth}}
  \subcaptionbox{\label{fig:porous_displacement}}{%
    \includegraphics[width=.37\textwidth,height=.37\textwidth]{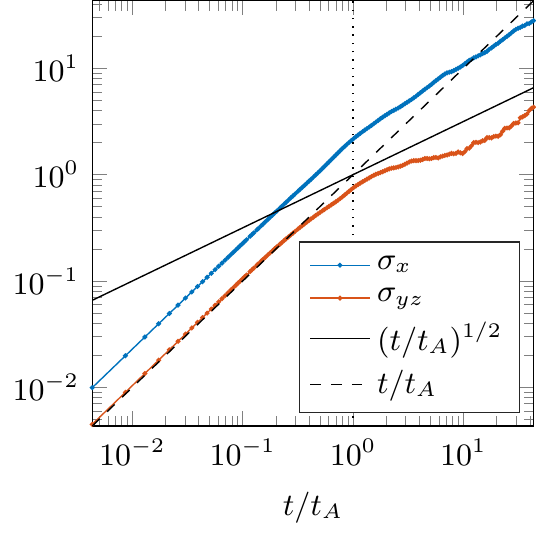}}
  \caption{(%
    \subref{fig:porous_example}) Streamlines around a periodic system
    with 131 oblate particles. %
    (\subref{fig:porous_displacement}) Mean-squared displacement in
    longitudinal ($\sigma_x$) and transversal ($\sigma_{yz}$)
    directions. The transition from linear to square root (Fickian)
    scaling appears to begin around $t=t_A$.}
  \label{fig:porous}
\end{figure}

Once we have the solution $\v q$ on the particles, it easy to evaluate
the velocity field anywhere in the domain. We can then place a
Lagrangian point randomly in the domain and let it be advected by the
flow, using a high-order ODE solver (MATLAB's \texttt{ode45}). Doing
this for a few particles as they pass through the domain allows us to
draw the streamlines of Figure \ref{fig:porous_example}. Doing this
for many ($\sim100$) particles over a long time interval ($t \in
[0,10^3]$), it is possible to extract statistics about the flow. We
shall here perform an analysis that closely follows that by de Anna et
al. \cite{DeAnna2013} for a two-dimensional medium with circular
grains. Our intention here is not to draw any new conclusions about
this particular case, only to show that our method makes the analysis
possible and that the results correspond to those in
\cite{DeAnna2013}.
 
The hydrodynamic particle dispersion can be measured through the
particle displacement $\Delta x_i(t) = x_i(t) - x_i(0)$.  Computing
the variance of the displacement%
\footnote{Brackets denote average over all Lagrangian particles.}%
, $\sigma_{x_i}^2(t) = \left\langle[\Delta x_i(t)-\langle\Delta
  x_i(t)\rangle]^2\right\rangle$, we can study the longitudinal
displacement in $\sigma_x^2$ and the transversal displacement in
$\sigma_{yz}^2 := \sigma_y^2 + \sigma_z^2$. These are shown in Figure
\ref{fig:porous_displacement}, with time normalized by $t_A =
\overline{d}/\overline{v}$, the mean particle diameter over the mean
velocity. Both $\sigma_x$ and $\sigma_{yz}$ can be seen to have an
initial linear growth, until entering a transition regime beginning at
$t=t_A$, after which the growth approaches a rate proportional to
$t^{1/2}$ (this is known as the Fickian regime). This is consistent
with \cite{DeAnna2013}.

\begin{figure}[htbp]
  \centering
  \begin{subfigure}[b]{0.49\textwidth}
    \centering
    \includegraphics[width=\textwidth,height=.65\textwidth]{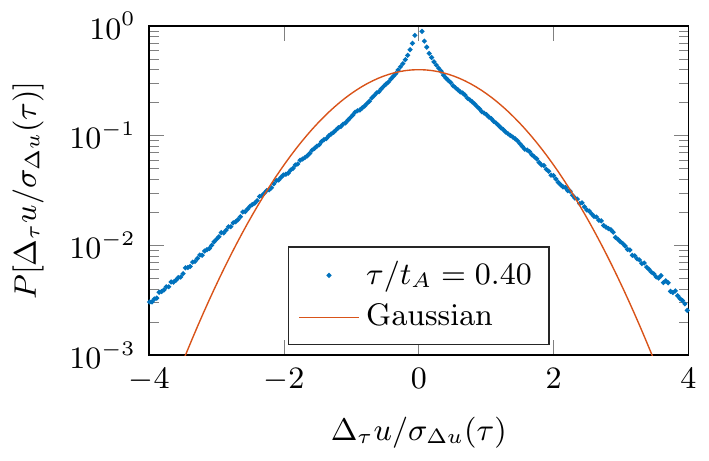}
  \end{subfigure}
  \begin{subfigure}[b]{0.49\textwidth}
    \centering
    \includegraphics[width=\textwidth,height=.65\textwidth]{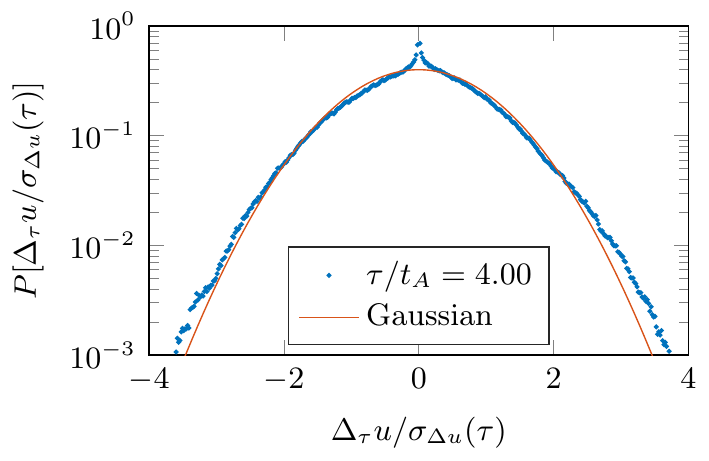}
  \end{subfigure}
  \caption{Probability density function of the normalized longitudinal Lagrangian velocity
    increments $\Delta_\tau u / \sigma_{\Delta u}(\tau)$ when $\tau < t_A$ and $\tau >
    t_A$. Clearly, a transition occurs around $\tau=t_A$.}
  \label{fig:porous_stats}
\end{figure}

For a given time lag $\tau$, the longitudinal Lagrangian velocity increment is defined as
$\Delta_\tau u = u(t+\tau) - u(t)$.  The probability density function of the velocity
increment, normalized by its standard deviation, is shown in Figure \ref{fig:porous_stats}
for two different time lags. For large time lags ($\tau \gg t_A$) the distribution
approaches the Gaussian distribution, while for small time lags ($\tau < t_A$) the
distribution is similar to that of the correlated continuous time random walk (CCTRW)
model in de Anna \cite{DeAnna2013}.

\section{Conclusions}
\label{sec:conclusions}

We have in this paper presented a complete boundary integral framework
for simulating periodic Stokes flow around spheroidal bodies. It is
based on representing the flow using the Stokes double layer
potential, which results in a well-conditioned system that converges
rapidly when solving using GMRES. Singular and nearly singular
integrals are computed using QBX, which we have adapted for the double
layer potential. This allows us to evaluate the velocity to high
accuracy everywhere in domain, such that we can have accurate
quadrature also for nearly touching geometries. By exploiting the
symmetries of the spheroidal bodies, we have been able to develop a
precomputation scheme that minimizes storage requirements, while
allowing us to rapidly compute the potential using QBX. By integrating
the precomputed QBX with a fast and spectrally accurate Ewald
summation method, we get a fast method which scales favorably as we go
to systems with many particles.

The results presented in section \ref{sec:results} show that our
method is able to accurately solve reference problems available in the
literature. The example involving a porous media model also indicates
that our method could be useful when studying such models, with both
spherical and non-spherical grains. A natural application extension of
our method would be to use it also for moving rigid particles, for
example by studying sedimentation in particle suspensions, as was done
in \cite{AfKlinteberg2014a}. Then particles can come arbitrary close
to each other, so close that lubrication forces start to
dominate. This is not an issue for our quadrature, but may cause the
double layer density to be sharply peaked, as we saw in the example in
section \ref{sec:triangle-spheres}. To deal with such cases, one could
either (i) introduce an adaptive surface quadrature, that is able to
resolve the density peak using a limited number of surface nodes, or
(ii) couple the method to a lubrication model that removes the need
for the double layer potential to resolve the strong interactions
\cite{Lefebvre-Lepot2015, Mammoli2006, Nasseri2000a, Sangani1994a}.

\section*{Acknowledgments}

This work has been supported by the Swedish Research Council under
grant no. \mbox{2011-3178}.

\appendix
\begin{appendices}

\section{Precomputing QBX on a spheroidal grid}
\label{sec:qbx-spheroidal}

In this appendix we develop the details of how to accelerate QBX for
the double layer potential on a spheroidal grid by exploiting
symmetries of the grid. To derive the expressions that we need, we
first consider acceleration on a spheroidal grid for the Laplace
potentials.

\subsection{Symmetries in precomputed layer potentials}
\label{sec:general-symmetries}

We begin by considering a general, scalar layer potential from a radially symmetric kernel
$K$, integrated over a surface $\Gamma$ together with a scalar density $q$,
\begin{align}
  u(\vx) = \int_\Gamma K(\vx-\vy)q(\vy)\d S_{\vy} .
\end{align}
The surface is discretized using a spheroidal grid, as defined in section
\ref{sec:spheroidal-grid}, which gives us $N=\ntheta\nphi$ points $\vx_1, \dots, \vx_N$ on
the surface. At each point we have a discrete density $q_i = q(\vx_i)$. All density values
are stored in the vector $\v Q \in \mathbb{R}^{N}$,
\begin{align}
  \v Q = 
  \begin{pmatrix}
    q(\vx_1) \\
    \\
    \\
    \vdots \\
    \\
    \\
    q(\vx_N)
  \end{pmatrix}
  = 
  \begin{pmatrix}
    q(\theta_1,\varphi_1) \\
    q(\theta_2,\varphi_1) \\
    \vdots \\
    q(\theta_\ntheta,\varphi_1) \\
    q(\theta_1,\varphi_2) \\
    \vdots \\
    q(\theta_\ntheta,\varphi_\nphi) \\
  \end{pmatrix} .
\end{align}
The discrete approximation of the layer potential is then
\begin{align}
  u^h(\vx) = \sum_{i=1}^N K(\vx-\vx_i)q_iw_i.
\end{align}
When solving a boundary integral equation, we typically need to compute the potential at
all the points $\v x_i$ on the surface,
\begin{align}
  u^h(\vx_i) = \sum_{j=1}^N K(\vx_i-\vx_j)q_jw_j, \quad i=1,\dots,N,
\end{align}
which we can write more compactly as
\begin{align}
  \v U = A \v Q,
\end{align}
where $A \in \mathbb R^{N \times N}$ is a (typically) full matrix and $\v U \in \mathbb
R^N$ is the layer potential at all the points on the surface, $U_i = u(\v x_i)$.

\subsubsection{Rotational symmetry}
Now, since the spheroid is rotationally symmetric and we use the trapezoidal rule with
uniform spacing in the azimuthal direction, the grid will look identical from the point of
view of all points which lie on the same latitude ($\theta$). If we denote the rows of $A$
by $R_i$ , $i=1,\dots,N$, then this must mean that rows $i$ and $i+\ntheta$ contains the
same information, and differ only by a permutation. This can be viewed as $R_i$ and
$R_{i+\ntheta}$ being stencils that cover the spheroidal grid in the
$\theta$-$\varphi$-plane. Having their target point at the same $\theta$-coordinate they
then have identical coefficients, but they differ by a permutation due to the periodic
wrap-around at $\phi=2\pi$.

We introduce the permutation $\gv{\tau}_n(\v Q)$ \footnote{$\textrm{mod}(a,b)$ refers to
  the operation $a$ \emph{modulo} $b$, which returns the remainder of the division of $a$
  by $b$, s.t. $0 \le \textrm{mod}(a,b) < b$.},
\begin{align}
  (\gv{\tau}_n(\v Q))_i = Q_{\sigma_n(i)}, \quad \sigma_n(i) = (i+n-1 \mod N)+1,
\end{align}
such that
\begin{align}
  \gv{\tau}_n(\v Q) = 
  \begin{pmatrix}
    Q_{1+n} \\
    \vdots \\
    Q_{N} \\
    Q_1 \\
    \vdots \\
    Q_n
  \end{pmatrix}.
\end{align}
We represent this permutation using a (sparse) $N\times N$ matrix $P_n$,
\begin{align}
  \gv{\tau}_n(\v Q) = P_n \v Q .
\end{align}
For $1 \le i \le \ntheta$ we then have that
\begin{align}
  u(\vx_i) &= R_i\vq, \\
  u(\vx_{i+\ntheta}) &= R_iP_\ntheta\vq, \\
  u(\vx_{i+2\ntheta}) &= R_iP_\ntheta^2\vq, \\
  &\dots ,
\end{align}
which we can summarize as
\begin{align}
    R_{i+\alpha\ntheta} &= R_iP_\ntheta^\alpha, \quad i \in [1,\ntheta] .
\end{align}
So we only need to store the first $\ntheta$ rows of $A$, as the remaining ones are given
by successive permutations of that set. This reduces our memory use from $N \times N$ to
$\ntheta \times N$ .

\subsubsection{Mirror symmetry}
For the rotational symmetry we could show that we only need to store coefficients related
to the first $\ntheta$ points in order to make computations for all $\ntheta \times \nphi$
points. An additional reduction can be made by observing that the spheroidal grid also has
a mirror symmetry around its equator. This means that we can use coefficients related to
target points on the northern hemisphere to make computations for target points on the
southern hemisphere. Hence, we only need to store the first $\ntheta/2$ rows of the matrix
$A$ in memory to be able to compute $\v U$ from $\v Q$.

Mathematically this is done for the current ordering of points by
introducing a mirroring permutation matrix $F_\ntheta$ such that
\begin{align}
  F_\ntheta\v Q =
  \begin{pmatrix}
    Q_\ntheta \\
    \vdots \\
    Q_1 \\
    Q_N \\
    \vdots \\
    Q_{\ntheta+1}
  \end{pmatrix} .
\end{align}
We can then recover rows $\ntheta/2+1,\dots,\ntheta$ by using rows $1,\dots,\ntheta/2$,
\begin{align}
  R_i = R_{\ntheta-i+1}F_\ntheta, \quad i \in [\ntheta/2+1,\ntheta],
\end{align}
and from there recover all the rows of the matrix using the rotational symmetry. For
notational simplicity we here assume $\ntheta$ to be even, extending to odd $\ntheta$ is
trivial.

\subsection{Precomputing QBX}

Now we have established the symmetries available when evaluating a
layer potential on a spheroidal grid. We will now see how those can be
used in conjunction with QBX, to create a fast method.

\subsubsection{Laplace single layer potential}
\label{sec:laplace-single-layer}

We begin by considering the Laplace single layer potential, having kernel $K(\v r) =
|r|^{-1}$, as the analysis is simplified by all the involved quantities being
scalar. Using the Laplace expansion \eqref{eq:laplace_expansion} about an expansion center
$\v c$, the potential
\begin{align}
  u(\v x) = \int_\Gamma \frac{q(\v y)}{|\v x-\v y|}\dSy,
\end{align}
can be written as
\begin{align}
  u(\v x) = \sum_{l=0}^\infty
   \sum_{m=-l}^l
  r_x^l Y_l^{-m} (\theta_x,\varphi_x) z_{lm}(\v c),
\end{align}
where the expansion coefficients are computed as
\begin{align}
  z_{lm}(\v c) = \frac{4\pi}{2l+1} \int_\Gamma \frac{q(\v y)}{r_y^{l+1}}
  Y_l^{m}(\theta_y,\varphi_y) \dSy .
\end{align}
When this is done discretely on a spheroidal grid, there is a series of steps
involved. The first step is to create a fine grid by increasing the number of grid points
by a factor $\kappa$ in each direction, so that we get a grid with points $\v{\widetilde
  x}_i$, $i=1,\dots,\kappa^2N$. The density is upsampled to the fine grid by a suitable
interpolation scheme, which is represented as the matrix operation
\begin{align}
  \widetilde{\v Q} = U \v Q, \quad U \in \mathbb{R}^{\kappa^2N\times N}.
\end{align}
In the second step, we use the upsampled density to compute the local expansion
coefficients up to order $p$ at expansions points $\v c_i$ (located at a distance $r$ in
the normal direction from $\vx_i$),
\begin{align}
  z_{lm}(\v c_i) = \sum_{n=1}^{\kappa^2N} \widetilde q_n \widetilde w_n 
  \frac{1}{r_n^{l+1}} Y_l^{m} (\theta_n,\varphi_n),
  \label{eq:def_sgl_mom}
\end{align}
where $\widetilde w_n$ is the quadrature weight at $\widetilde \vx_n$ and
$(r_n,\theta_n,\varphi_n) = \widetilde \vx_n - \v c$. This gives us for each point $\v
c_i$ a coefficient vector $\v z(\v c_i) \in \mathbb{C}^{{\Nmom}}$, 
\begin{align}
  \v z(\v c_i) =
  \begin{pmatrix}
    z_{0,0}(\v c_i) \\
    z_{1,0}(\v c_i) \\
    z_{1,1}(\v c_i) \\
    \vdots \\
    z_{p,p}(\v c_i)
  \end{pmatrix} ,  
\end{align}
where ${\Nmom} = (p^2+3p+2)/2$ is the total number of coefficients at each expansion
center that we need to store\footnote{As mentioned in section \ref{sec:qbx-stresslet}, the
  coefficients for negative $m$ need not be stored, since $z_{l,-m}=z_{lm}^*$ for a real
  density.}. In matrix notation we write this as
\begin{align}
  \v z(\v c_i) = \mathcal{S}_i \widetilde{\v Q} =  \mathcal{S}_i U \v Q = M_i \v Q,
\label{eq:sgl_den2mom}
\end{align}
where
\begin{align}
  \mathcal{S}_i &\in \mathbb{C}^{{\Nmom} \times \kappa^2 N}, \\
  M_i &\in \mathbb{C}^{{\Nmom} \times N} .
\end{align}
As the third step, we let the vector $\v s(\v x_i)\in\mathbb{C}^{{\Nmom}}$ contain the
coefficients for evaluating the expansion at $\v c_i$,
\begin{align}
  \v s(\v s_i) &=
  \begin{pmatrix}
    s_{0,0}(\v x_i) \\
    s_{1,0}(\v x_i) \\
    s_{1,1}(\v x_i) \\
    \vdots \\
    s_{p,p}(\v x_i)
  \end{pmatrix} ,&  
  s_{lm}(\v x_i) &= 
  \begin{cases}
    \frac{1}{2} r_x^l Y_l^{m} (\theta_x,\varphi_x), & \text{ if } m=0, \\
    r_x^l Y_l^{m} (\theta_x,\varphi_x) & \text{ if } m>0
  \end{cases},
  \label{eq:sg_mom2pot_coeff}
\end{align}
such that, using $(Y_l^m)^*=Y_l^{-m}$, we can write
\begin{align}
  u(\v x_i) = \sum_{l=0}^p \sum_{m=0}^l 
  ( s_{lm}(\v x_i)^* z_{lm}(\v c) + s_{lm}(\v x_i) z_{lm}(\v c)^*) .
\end{align}
Defining the conjugate product \newcommand{\conjprod}[2]{\langle#1,#2\rangle}
\begin{align}
  \conjprod{\vx}{\vy} = (\vx^*)^T\vy + \vx^T\vy^*,
\end{align}
we can write this as
\begin{align}
  u(\vx_i) = \conjprod{\v s(\vx_i)}{\v z(\v c_i)}.
  \label{eq:mom2pot}
\end{align}
Putting \eqref{eq:sgl_den2mom} and \eqref{eq:mom2pot} together, we can conclude that the
matrix $R_i$, such that $u(\v x_i) = R_i\v q$, is given by
\begin{align}
  R_i = \v s(\vx_i)^{*,T} M_i + \v s(\vx_i)^T M_i^*, \quad i=1,\dots,N .
\end{align}
This $R_i$ has the symmetries discussed in section \ref{sec:general-symmetries}. From
rotation and mirroring when then have
\begin{align}
  R_i &= R_{\ntheta-i+1}F_\ntheta, &i &\in [\ntheta/2+1,\ntheta], \\
  R_{i+\alpha\ntheta} &= R_iP_\ntheta^\alpha,& i &\in [1,\ntheta],
\end{align}
such that it is enough to directly compute $R_i,\dots,R_{\ntheta/2}$, as the remaining $R_i$
can be recovered through these operations.

The symmetries that apply to the rows $R_i$ also apply to the matrices $M_i$, such that we
can compute the coefficients at all expansion centers using $M_i$ for the first $\ntheta$
expansion centers. Here, however, it is not enough to permute the right hand side vector;
we also need to account for the phase of the factor $e^{im\varphi}$ in the spherical
harmonic. An expansion center $\v c_{i+\ntheta}$ will ''see'' all the grid points in the
same angles $\varphi$ as an expansion center $\v c_i$ does, plus a rotation of
$\Delta\varphi$. To account for this rotation, we introduce a phase matrix $E \in
\mathbb{C}^{{\Nmom}\times {\Nmom}}$ with diagonal elements
\begin{align}
  E_{jj} = e^{im\Delta\varphi}.
\end{align}
For $1 \le i \le \ntheta$ this lets us write
\begin{align}
  \v z(\v c_i) &= M_i \v Q, \\
  \v z(\v c_{i+\ntheta}) &= E M_i P_\ntheta \v Q, \\
  \v z(\v c_{i+2\ntheta}) &= E^2 M_i P_\ntheta^2 \v Q, \\
  &\dots,
\end{align}
and
\begin{align}
  \v s(\vx_{i+\ntheta}) &= \v s(\vx_i) E \\
  \v s(\vx_{i+2\ntheta}) &= \v s(\vx_i) E^2 \\
  &\dots,
\end{align}
which can be summarized as 
\begin{align}
    M_{i+\alpha\ntheta} &= E^\alpha M_i P_\ntheta^\alpha, \\
    \v s(\vx_{i+\alpha\ntheta}) &= \v s(\vx_i) E^\alpha .
\end{align}
From this we can recover the original rotational symmetry of $R_i$,
\begin{align}
  R_{i+\alpha\ntheta} &= 
  \v s(\vx_{i+\alpha\ntheta})^* M_{i+\alpha\ntheta} + \v s(\vx_{i+\alpha\ntheta}) M_{i+\alpha\ntheta}^* \\
  &=
  \left(\v s(\vx_i) E^\alpha\right)^* E^\alpha M_i P_\ntheta^\alpha
  + \v s(\vx_i) E^\alpha  \left(E^\alpha M_i P_\ntheta^\alpha\right)^* \\
  &=
  \v s(\vx_i)^* M_i P_\ntheta^\alpha + \v s(\vx_i) M_i^* P_\ntheta^\alpha \\
  &=
  R_iP_\ntheta^\alpha ,
\end{align}
since $E^*E=EE^*=I$.

To also make use of the mirror symmetry for $M_i$, we have to keep in mind that the
mirroring operation makes an expansion center $\v c_{\ntheta/2+i}$ ''see'' the other grid
points in a coordinate system which has been rotated $180^\circ$ about the $x$-axis,
compared to what $\v c_i$ sees. This amounts to $\theta$ and $\varphi$ in
\eqref{eq:def_sgl_mom} having counterparts $\theta'=-\theta$ and $\varphi'=-\varphi$ in
the mirrored system. Inserting this into the relevant terms in \eqref{eq:def_sgl_mom}, we
see that
\begin{align}
  Y_l^m(\theta',\varphi') = 
  \widetilde P_l^m(\cos\theta')e^{im\varphi'} =
  \widetilde P_l^m(-\cos\theta)e^{-im\varphi} .
\end{align}
Changing sign in the exponential is equal to taking the complex
conjugate of the entire expression. To handle the change of sign in
the argument of the associated Legendre polynomial, we apply the
parity relation available in Arfken \& Weber \cite[p. 776]{Arfken2005},
\begin{align}
  P_l^m(-x) = (-1)^{l+m}P_l^m(x) .
  \label{eq:assleg_parity}
\end{align}
\newcommand{\mirrtheta}{{\Theta}}
Defining the diagonal matrix $\mirrtheta \in \mathbb{R}^{{\Nmom} \times {\Nmom}}$ such
that
\begin{align}
  \mirrtheta_{jj} = (-1)^{l+m}, 
\end{align}
where $l=l(j)$ and $m=m(j)$ as $j=1,\dots,{\Nmom}$, we finally arrive at
the relation
\begin{align}
    M_i = \mirrtheta M_{\ntheta-i+1}^*F_\ntheta, \quad i \in [\ntheta/2+1,\ntheta],,
\end{align}
which allows us to compute all expansion coefficients for all expansion
points using precomputed results for the first $\ntheta/2$ points.

\subsubsection{Dipole potential}
\label{sec:dipole-potential}

We will now extend the above results for the Laplace single layer potential to the dipole
potential, defined in \eqref{eq:ldbl_def} for a vector density $\v q$ as
\begin{align}
  \ldbl{\Gamma,\v q}{\v x} := 
  \int_\Gamma \v q(\v y) \cdot\nabla_{\v y}\frac{1}{|\v x - \v y|} \dSy 
  .
\end{align}
A discrete quadrature by expansion for this potential is formed from the expressions
\eqref{eq:ldbl_expa} and \eqref{eq:ldbl_exp_coeff}. We store the pointwise vector
densities in linear form in $\v Q \in \mathbb R^{3N}$,
\begin{align}
  \v Q = 
  \begin{pmatrix}
    \v Q_1 \\ \v Q_2 \\ \v Q_3
  \end{pmatrix},
  \quad\text{where}\quad
  \v Q_j =
  \begin{pmatrix}
    q_j(\v x_1) \\
    \vdots \\
    q_j(\v x_N)
  \end{pmatrix} .
  \label{eq:linear_vector_density}
\end{align}
Just as for the single layer potential, the expansion coefficients are computed using an
upsampled density $\v{\widetilde Q}$, which here is computed as
\begin{align}
  \v{\widetilde Q} = (I_3\otimes U)\v Q,
\end{align}
where $\otimes$ is the Kronecker product. The expansion coefficients
are here computed as
\begin{align}
  z_{lm}(\v c) = \sum_{n=1}^{\kappa^2N} \widetilde w_n
  \v{\widetilde q}(\v{\widetilde x}_n) \cdot
  \nabla_{\widetilde \vx} \frac{1}{r^{l+1}} \widetilde
  Y_l^m(\theta_n, \varphi_n) .
\end{align}
This can be represented in matrix form as
\begin{align}
  \v z(\v c_i) =\mathcal{S}_i (I_3\otimes U) \v Q = M_i\v Q,
  \label{eq:dipole_den2mom}
\end{align}
where
\begin{align}
  \v z(\v c_i) &\in \mathbb C^{N_p}, \\
  M_i &\in \mathbb C^{N_p \times 3N}.
\end{align}
The difference from the single layer case is that the matrix $\mathcal{S}_i \in \mathbb
C^{\Nmom \times 3\kappa^2N}$ here includes the gradient of the spherical harmonic on the
surface. Evaluation of the expansion is identical to the single layer case
\eqref{eq:mom2pot},
\begin{align}
  u(\vx_i) = \conjprod{\v s(\vx_i)}{\v z(\v c_i)},
  \label{eq:dbl_mom2pot}
\end{align}
with $\v s(\v x_i)$ defined as in \eqref{eq:sg_mom2pot_coeff}, and $R_i$ given by
\begin{align}
  R_i = \v s(\vx_i)^{*,T} M_i + \v s(\vx_i)^T M_i^*, \quad i=1,\dots,N .
\end{align}

Due to the density now being vector-valued, there is difference between the single layer
potential and the dipole potential in how we can make use of the grid symmetries. To reuse
$M_i$ of the first $\ntheta$ points by rotational symmetry we must now not only take the
permutation $P_\ntheta$ into account, but also the rotation of the frame of reference,
since $\v Q$ now represents a quantity which is pointwise vector-valued. Introducing the
rotation matrix $T_z$,
\begin{align}
  T_z(\alpha) = 
  \begin{pmatrix}
    \cos(\alpha) & -\sin(\alpha) & 0 \\
    \sin(\alpha) & \cos(\alpha) & 0 \\
    0 & 0 & 1
  \end{pmatrix},
\end{align}
the rotation of $M_i$ can be formulated as
\begin{align}
  M_{i+\alpha\ntheta} = E^\alpha M_i \left( T_z(-\Delta\varphi) \otimes
    P_\ntheta \right)^\alpha .
  \label{eq:M_rotsym_dip}
\end{align}
Similarly, for the mirror symmetry we must take the $180^\circ$
rotation about the $x$-axis into account. Representing this by the
rotation matrix
\begin{align}
  T_x = 
  \begin{pmatrix}
    1 & 0 & 0 \\
    0 &-1 & 0 \\
    0 & 0 &-1
  \end{pmatrix},
\end{align}
we arrive at the expression
\begin{align}
  M_i = \mirrtheta M_{\ntheta-i+1}^* (T_x \otimes F_\ntheta), \quad i \in
  [\ntheta/2+1,\ntheta],
  \label{eq:M_mirsym_dip}
\end{align}
for the mirror symmetry. These results also carry over to the symmetry relations for
$R_i$,
\begin{align}
  R_i &= R_{\ntheta-i+1} (T_x \otimes F_\ntheta), & i &\in [\ntheta/2+1,\ntheta] ,
  \label{eq:R_mirsym_dip}
  \\
  R_{i+\alpha\ntheta} &= R_i \left( T_z(-\Delta\varphi) \otimes
    P_\ntheta \right)^\alpha, & i &\in [1,\ntheta] .
  \label{eq:R_rotsym_dip}
\end{align}

\subsubsection{Gradient of expansion}
\label{sec:gradient-expansion}

To compute the gradient of a (single or double layer) potential that is locally
represented by a local expansion at $\v c$, we differentiate \eqref{eq:mom2pot} or
\eqref{eq:dbl_mom2pot} with respect to $\vx$,
\begin{align}
  \nabla u(\vx) =  \conjprod{\v g(\vx)}{\v z(\v c)},
\end{align}
where $\v g(\v x) \in \mathbb{C}^{3 \times {\Nmom}}$ has elements
\begin{align}
  g_{ij}(\vx) = \pd{}{x_i} s_j(\vx).
  \label{eq:def_g}
\end{align}
We can also write
\begin{align}
  \pd{}{x_i} u(\vx) =  \conjprod{\v g_i(\vx)}{\v z(\v c)},
\end{align}
where $\v g_i(\v x) \in \mathbb C^\Nmom$,
\begin{align}
  \v g_i(\v x) = \pd{}{x_i} \v s_i(\v x).
  \label{eq:def_gi}
\end{align}
If we want to reuse $\v g(\v x_i)$ from the first $\ntheta$ points
for the remainder of the points, we have to take into account that the
output is vector-valued and has to be rotated into the correct frame
of reference. So, if we have the entire process condensed into the form
\begin{align}
  \nabla u(\vx_i) = R_i \v q, \quad 1 \le i \le \ntheta,
\end{align}
then
\begin{align}
  \nabla u(\vx_{i+\alpha\ntheta}) = T_z^\alpha(\Delta\varphi)R_iP_\ntheta^\alpha \v q, \quad
  1 \le \alpha < \nphi.
  \label{eq:R_rotsym_grad}
\end{align}
For the mirror symmetry we get
\begin{align}
  \nabla u(\vx_{\ntheta/2+i}) = T_xR_{\ntheta/2-i+1}F_\ntheta, \quad 1 \le i \le \ntheta/2 .
  \label{eq:R_mirsym_grad}
\end{align}

\subsubsection{Stokes double layer potential}
\label{sec:stokes-double-layer}
We now turn to our potential of interest,
the Stokes double layer potential \eqref{eq:dbl_def},
\begin{align}
  \dbli{\Gamma,\v q}{\v x} = \int_{\Gamma} \stresslet_{ijk}(\v x,\v y) q_j(\v y)
  \nhat_k(\v y) \dSy.
\end{align}
Our goal is now to derive the matrix representations necessary for precomputing QBX for
this potential on a spheroidal grid. We have from \eqref{eq:stresslet_decomp} that the we
can represent the double layer potential as a combination of dipole potentials,
\begin{align}
  \begin{split}
    \dbli{\Gamma,\v q}{\v x} =
    & \left( x_j\pd{}{x_i}-\delta_{ij} \right) \ldbl{\Gamma,q_j
      \v\nhat+\nhat_j \v q}{\v x}
    \\
    & - \pd{}{x_i}\ldbl{\Gamma,y_kq_k\v\nhat+y_k\nhat_k\v q}{\v x}    
  \end{split} ,
  \label{eq:stresslet_decomp_repeat}
\end{align}
which means that we can derive the expressions needed for the double layer potential using
the results derived so far in sections
\ref{sec:laplace-single-layer}--~\ref{sec:gradient-expansion}. Taking
\eqref{eq:stresslet_decomp_repeat} apart, we get
\begin{align}
  \dbli{\Gamma,\v q}{\v x} = \left( x_j\pd{}{x_i}-\delta_{ij} \right) u^j(\v x) -
  \pd{}{x_i} u^4(\v x),
\end{align}
where
\begin{align}
  u^j(\v x) &= \ldbl{\Gamma,q_j \v\nhat+\nhat_j \v q}{\v x}, & j&=1,\dots,3, \\
  u^4(\v x) &= \ldbl{\Gamma,y_kq_k\v\nhat+y_k\nhat_k\v q}{\v x} . &&    
\end{align}
When computing these components on a spheroidal grid, we can improve the accuracy by
taking into account that the geometrical quantities $\v\nhat$ and $\vy$ are explicitly
known, such that we only need to upsample the density $\vq$. Using
\eqref{eq:dipole_den2mom}, the expansion coefficients at a center $\v c$ for the
components $u^1$--$u^4$ are then given by
\begin{align}
  \v z^j(\v c) &= \mathcal{S}(\v c) \left [
    \widetilde{\v N}^T
    U \v Q_j
    + \left(
      I_3 \otimes (\text{diag}(\v{\widetilde N}_j) U)
    \right) \v Q
    \right ] ,& j &= 1,2,3,
    \label{eq:dbl_zj}
  \\
  \v z^4(\v c) &= \mathcal{S}_i(\v c) \left[
    \v{\widetilde  N}^T
    \v{\widetilde X}
    (I_3 \otimes U)
    + (I_3 \otimes \v{\widetilde P} U)
  \right ] \v Q, & &
  \label{eq:dbl_z4}
\end{align}
where the geometrical quantities on the fine grid are given by
\begin{align}
  \v{\widetilde X}_j &=
  \begin{pmatrix}
    (\v{\tilde x}_{1})_j \\
    \vdots \\
    (\v{\tilde x}_{\kappa^2N})_j
  \end{pmatrix}
  ,&
  \v{\widetilde  N}_j &=
  \begin{pmatrix}
    (\v{\tilde \nhat}_{1})_j \\
    \vdots \\
    (\v{\tilde \nhat}_{\kappa^2N})_j
  \end{pmatrix}
\end{align}
and
\begin{align}
  \v{\widetilde N} &=
  \begin{pmatrix}
    \text{diag}(\v{\widetilde  N}_1) &
    \text{diag}(\v{\widetilde  N}_2) &
    \text{diag}(\v{\widetilde  N}_3) 
  \end{pmatrix}, \\
  \v{\widetilde X} &=  
  \begin{pmatrix}
    \text{diag}(\v{\widetilde X}_1) &
    \text{diag}(\v{\widetilde X}_2) & 
    \text{diag}(\v{\widetilde X}_3) 
  \end{pmatrix}, \\
  \v{\widetilde P} &= \text{diag}\left( 
    \v{\tilde x}_1 \cdot \v{\tilde\nhat}_{1}, \dots,
    \v{\tilde x}_{\kappa^2N} \cdot \v{\tilde\nhat}_{\kappa^2N} 
  \right) .
\end{align}
From the coefficient vectors $\v z^1$--$\v z^4$ at $\v c$, the double layer potential at a
point $\v x$ can be computed as
\begin{align}
  \dbli{\Gamma,\v q}{\v x} &= \conjprod{x_j\v g_i(\vx)-\delta_{ij}\v s(\vx)}{\v z_j(\v c)}
  -\conjprod{\v g_i(\vx)}{\v z_4(\v c)},
  \label{eq:dbl_eval}
\end{align}
where the evaluation vectors $\v s(\v x)$ and $\v g_i(\v x)$ are as defined in
\eqref{eq:sg_mom2pot_coeff} and \eqref{eq:def_gi}.

For the points $\v x_i$ and centers $\v c_i$ related to the spheroidal grid, the
operations \eqref{eq:dbl_zj},~\eqref{eq:dbl_z4} and \eqref{eq:dbl_eval} can be represented 
using a set of transfer matrices $M_i^{1-4}\in\mathbb{C}^{{\Nmom} \times 3N}$ and
$R_i\in\mathbb{R}^{3 \times 3N}$, such that
\begin{align}
  \v z^j(\v c_i) &= M_i^jQ, \quad j=1,2,3,4,\\
  \dbl{\Gamma,\v q}{\vx_i} &= R_iQ,
\end{align}
where $i=1,\dots,N$. We can, just as for the single layer and dipole potentials, compute
the matrices $R_i$ and $M_i^j$ for all points by applying rotation and mirror operations
to the matrices from the first $\ntheta/2$ points. For $R_i$, we can combine
\eqref{eq:R_mirsym_dip},~\eqref{eq:R_rotsym_dip},~\eqref{eq:R_rotsym_grad} and
\eqref{eq:R_mirsym_grad} directly to get
\begin{align}
  R_i &= T_x R_{\ntheta-i+1} (T_x \otimes F_\ntheta), & i &\in [\ntheta/2+1,\ntheta], \\
  R_{i+\alpha\ntheta} &= T_z^\alpha(\Delta\varphi) R_i \left( T_z(-\Delta\varphi) \otimes
    P_{\ntheta} \right)^\alpha, & i &\in [1, \ntheta] .
\end{align}
For the matrices $M_i^j$, we need to take into account that $M_i^1$--$M_i^3$ are
directly related to the coordinate directions $x$,$y$ and $z$. This means that when
applying the symmetry operations to these matrices, we have to take this into account. For
the mirror symmetry we then get
\begin{align}
  \begin{split}
    M_i^1 &= \mirrtheta M_{\ntheta-i+1}^{1,*} (T_x \otimes F_\ntheta),\\
    M_i^2 &= -\mirrtheta M_{\ntheta-i+1}^{2,*} (T_x \otimes F_\ntheta),\\
    M_i^3 &= -\mirrtheta M_{\ntheta-i+1}^{3,*} (T_x \otimes F_\ntheta),\\
    M_i^4 &= \mirrtheta M_{\ntheta-i+1}^{4,*} (T_x \otimes F_\ntheta),\\
  \end{split}
  \quad i \in [\ntheta/2+1,\ntheta],
\end{align}
while for the rotational symmetry
\begin{align}
  \begin{split}
    M_{i+\alpha\ntheta}^{1} &= E^\alpha \left( \cos(\alpha\Delta\varphi) M_{i}^{1} -
      \sin(\alpha\Delta\varphi) M_{i}^{2} \right) \left( T_z(-\Delta\varphi) \otimes
      P_\ntheta
    \right)^\alpha, \\
    M_{i+\alpha\ntheta}^{2} &= E^\alpha \left( \sin(\alpha\Delta\varphi) M_{i}^{1} +
      \cos(\alpha\Delta\varphi) M_{i}^{2} \right) \left( T_z(-\Delta\varphi) \otimes
      P_\ntheta
    \right)^\alpha, \\
    M_{i+\alpha\ntheta}^{3} &= E^\alpha M_{i}^{3} \left( T_z(-\Delta\varphi) \otimes P_\ntheta
    \right)^\alpha, \\
    M_{i+\alpha\ntheta}^{4} &= E^\alpha M_{i}^{4} \left( T_z(-\Delta\varphi) \otimes P_\ntheta
    \right)^\alpha .
  \end{split}
\end{align}
With this we have all we need for using precomputed QBX on a spheroidal grid, while only
storing precomputed data for the first $\ntheta/2$ points.

\end{appendices}

\clearpage

\bibliographystyle{jabbrv_abbrv_doi}
\bibliography{../../Documents/Library/library}

\end{document}

%% file: defs.tex
\renewcommand{\v}[1]{\bm{#1}}
\newcommand{\gv}[1]{\bm{#1}} 

\renewcommand{\d}{\mathrm{d}}
\newcommand{\pd}[2]{\frac{\partial #1}{\partial #2}} 
\DeclareMathOperator{\erfc}{erfc}

\newcommand{\ordo}{\mathcal{O}}

\newcommand{\ntheta}{{n_{\theta}}}
\newcommand{\nphi}{{n_{\varphi}}}
\newcommand{\nhat}{n}
\newcommand{\dSy}{\d S_{\v y}}

\newcommand{\stokeslet}{S}
\newcommand{\stresslet}{T}
\newcommand{\rotlet}{R}
\newcommand{\dipole}{F}
\newcommand{\Npart}{M}
\newcommand{\gamalp}{\Gamma_{\alpha}}
\newcommand{\asup}{^{\alpha}}
\newcommand{\rspace}{^{(R)}}
\newcommand{\kspace}{^{(F)}}
\newcommand{\tol}{\epsilon}

\newcommand{\vx}{\v{x}}
\newcommand{\vy}{\v{y}}
\newcommand{\vq}{\v{q}}
\newcommand{\potential}[3]{\ensuremath{
    #1        
    \left[       
      #2         
    \right]
    \left(#3\right)
}}
\newcommand{\numapprox}[1]{#1^{h}}
\newcommand{\ldblsymb}{\mathcal{L}}
\newcommand{\ldbl}[2]{\potential{\ldblsymb}{#1}{#2}}
\newcommand{\dblsymb}{\ensuremath{\mathcal{D}}}
\newcommand{\dbl}[2]{\potential{\v\dblsymb}{#1}{#2}}
\newcommand{\dbli}[2]{\potential{\dblsymb_i}{#1}{#2}}
\newcommand{\dblh}[2]{\potential{\numapprox{\v\dblsymb}}{#1}{#2}}
\newcommand{\ucompsymb}{\mathcal{V}}
\newcommand{\ucomp}[1]{\ensuremath{ \v\ucompsymb\left(#1\right) }}
\newcommand{\ucompi}[1]{\ensuremath{ \ucompsymb_i\left(#1\right) }}
\newcommand{\ubg}{{\v{u}_\text{bg}}}
\newcommand{\mself}{R}
\newcommand{\mmom}{M}
\newcommand{\Nmom}{{N_p}}
\newcommand{\opint}{\operatorname{I}}
\newcommand{\opquad}{\operatorname{Q}}

%% file: acronyms.tex
\usepackage{acronym}

\newacro{BIE}[BIE]{boundary integral equation}

\newacro{FFT}[FFT]{fast Fourier transform}
\newacroindefinite{FFT}{an}{a}

\newacro{FMM}[FMM]{fast multipole method}
\newacroindefinite{FMM}{an}{a}

\newacro{GMRES}[GMRES]{the generalized minimum residual method \cite{Saad1986}}

\newacro{SE}[SE]{spectral Ewald}